\documentclass[11pt]{article}
\usepackage{graphicx,mathrsfs,amssymb,color}
\usepackage{amsmath,latexsym,amsbsy}
\usepackage{tikz}
\usepackage{pgfplots}
\usetikzlibrary{matrix}

\def\ds{\displaystyle}
\def\forall{\hbox{for all}~}

\def\ve{\varepsilon}
\def\n{\noindent}

\def\R{\mathbb{R}}

\def\vs{\vskip 2em}

\def\v{\vskip 1em}
\def\O{{\cal O}}

\def\bega{\begin{array}}
\def\enda{\end{array}}
\def\begi{\begin{itemize}}
\def\endi{\end{itemize}}

\def\bel{\begin{equation}\label}
\def\eeq{\end{equation}}
\def\sqr#1#2{\vbox{\hrule height .#2pt
\hbox{\vrule width .#2pt height #1pt \kern #1pt
\vrule width .#2pt}\hrule height .#2pt }}
\def\square{\sqr74}
\def\endproof{\hphantom{MM}\hfill\llap{$\square$}\goodbreak}

\newtheorem{theorem}{Theorem}[section]
\newtheorem{corollary}{Corollary}[section]
\newtheorem{lemma}{Lemma}[section]
\newtheorem{proposition}{Proposition}[section]
\newtheorem{remark}{Remark}[section]
\newtheorem{definition}{Definition}[section]

\usepackage{tikz}
\usepackage{placeins}
\usepackage{soul}

\usepackage{chngcntr}
\counterwithin{equation}{section}
\usepackage{enumitem}

\begin{document}
\title{\bf On the rate of convergence in superquadratic Hamilton–Jacobi equations with state constraints}\vs
\author{\it  Prerona Dutta$^{(1)}$, Khai T. Nguyen$^{(2)}$, and Son N.T. Tu$^{(3)}$\\
\\
{\small $^{(1)}$ Department of Mathematics, Xavier University of Louisiana}\\
		{\small $^{(2)}$  Department of Mathematics, North Carolina State University,}\\
		{\small $^{(3)}$  Department of Mathematics, Baylor University.}\\
		\\
		\quad\\
		{\small E-mails:  : pdutta@xula.edu,~~ khai@math.ncsu.edu,~~ son\_tu@baylor.edu}}
\maketitle

\begin{abstract} 

In this paper, we investigate the convergence rate in the vanishing viscosity limit for solutions to superquadratic Hamilton–Jacobi equations with state constraints. For every $p>2$, we establish the rate of convergence for nonnegative Lipschitz data vanishing on the boundary to be of order \( \mathcal{O}(\varepsilon^{1/2}) \)  and obtain an improved upper rate of order \( \mathcal{O}\big(\varepsilon^{{p\over 2(p-1)}}\big)\) for semiconcave data.
\v
		{\footnotesize
		\n {\bf Keywords.}
		Hamilton–Jacobi equations, problems with state constraints, viscosity solutions, rates of
convergence
		
		\medskip
		
		\n {\bf AMS Mathematics Subject Classification.}
	35B40 · 35D40 · 49J20 · 49L25 · 70H20		
		 
		}
	\end{abstract}

\section{Introduction}

\color{black}

\label{sec:1}
Throughout the paper, we consider an open, bounded domain $\Omega \subset \mathbb{R}^n$ with $\mathcal{C}^2$ boundary, and assume $p > 2$.
For every continuous function $f:\overline{\Omega}\to\R$,  we consider the following Hamilton-Jacobi  equation with state constraint:
\bel{HJ-S}
\begin{cases}
\lambda u(x)+H(x,D u(x))~=~0,\qquad x\in \Omega\\[2mm]
\lambda u(x)+H(x,D u(x))~\geq~0, \qquad x\in \partial\Omega
\end{cases}
\eeq
where  the Hamiltonian $H:\overline{\Omega}\times\R^n\to\R$ is defined by 
\bel{H}
H(x,\xi) = |\xi|^p - f(x),\qquad (x,\xi)\in \overline{\Omega}\times\R^n.
\eeq
The  Legendre transform $L$ of $H$ is computed explicitly as
\bel{L}
	L(x,v)~=~\sup_{\xi \in \mathbb{R}^n} \Big(\xi\cdot v - H(x,\xi)\Big) = c_p|v|^q + f(x), \qquad (x,v) \in \overline{\Omega}\times \mathbb{R}^n.
\eeq
where $q$ is the conjugate of $p$ such that $p^{-1}+q^{-1}=1$ and $c_p~=~ p^{-\frac{1}{p-1}} \left(1-\frac{1}{p}\right)$. Let  $\mathcal{AC}((-\infty,0];\overline{\Omega})$ denote  the set of absolutely continuous functions from $(-\infty,0]$ to $\overline{\Omega}$. It is well-known from \cite{soner_optimal_1986} that the unique viscosity solution of \eqref{HJ-S} is given by the representation formula
\begin{equation}\label{O-C}
	u(x)~=~\inf \left\lbrace 
	\int_{-\infty}^{0}e^{\lambda s} 
	L(\eta(s),\dot{\eta}(s))~ds: 
    \eta\in \mathcal{AC}((-\infty,0];\overline{\Omega}), \eta(0)=x
    \right\rbrace,
\end{equation}
We investigate the rate of convergence of $\ds\|u^{\ve}-u\|_{\mathcal{C}^0}$ as $\ve\to 0^+$, where $u^{\ve}$ is the value function of a stochastic optimal control problem subject to state constraints, solving  the second-order problem
\bel{HJ-ve}
\begin{cases}
\lambda u(x)+|D u(x)|^{p}-\ve\cdot \Delta u(x)~=~f(x),\qquad x\in \Omega \\[2mm]
\lambda u(x)+|D u(x)|^{p}-\ve\cdot \Delta u(x)~\geq~f(x),\qquad x\in \partial\Omega
\end{cases} ~.
\eeq

 For $1<p\leq 2$, equation \eqref{HJ-ve} is equivalent to
\begin{equation}\label{eq:SubQuadInf}
	\begin{cases}
	\begin{aligned}
	\lambda u(x) + |Du(x)|^p - \varepsilon \Delta u(x)  &~=~ 0,&& \qquad x\in \Omega\\[2mm]
	u(x) &~=~ \infty && \qquad x\in \partial\Omega
	\end{aligned}
	\end{cases}
\end{equation}
which is associated with a class of solutions called \emph{large solutions} in existing literature, namely \cite{porretta_local_2004,alessio_asymptotic_2006} and the references therein. For this case, the rate of convergence of solutions $u^\varepsilon\to u$ from \eqref{HJ-ve} to \eqref{HJ-S} was first studied in \cite{han_remarks_2022}. More precisely, it was shown that in addition to the rate of order $\mathcal{O}(\varepsilon^{1/2})$ for nonnegative continuous data vanishing on the boundary, improved rates of orders $\mathcal{O}(\varepsilon)$ and $\mathcal{O}(\varepsilon^{1/p})$ can be obtained for compactly supported and suitably regular semiconcave data respectively.

This paper establishes the  convergence rates for \( u^\varepsilon \to u \) in the case \( p > 2 \) where the behavior of the solution \( u^\varepsilon \) to \eqref{HJ-ve} near \( \partial\Omega \) is not explicitly known, unlike when \( p \leq 2 \). This lack of explicit boundary information makes it difficult to construct suitable barrier functions and the argument in  \cite{han_remarks_2022} can no longer be applied directly. Our first result is stated as follows:
\begin{theorem}\label{thm:NonnegativeZeroBoundaryIntro} For every $p>2$, assume that $\Omega\subset\R^n$ is an open and bounded  subset of $\R^n$ with $\mathcal{C}^2$ boundary and $f\in \mathbf{Lip}(\overline{\Omega})$. Then for every $\ve>0$, the following estimate holds:
\begin{align*}
\inf_{x\in\overline{\Omega}}\{u^\varepsilon(x) - u(x)\}~\geq~-{\underline{\mathbf{\Lambda}}\over\lambda} \cdot \varepsilon^{1/2}.  
\end{align*}
Additionally, if $f$ is nonnegative and vanishes on $\partial\Omega$, then
\[
\sup_{x\in\overline{\Omega}}\{u^\varepsilon(x) - u(x)\}~\leq~{\overline{\mathbf{\Lambda}}\over\lambda} \cdot \varepsilon^{1/2}. 
\]
\end{theorem}
 The proof of Theorem \ref{thm:NonnegativeZeroBoundaryIntro} is divided into Propositions \ref{prop:LowerBound} and Proposition \ref{prop:upper}, where the constants $\underline{\mathbf{\Lambda}} $ and  $\overline{\mathbf{\Lambda}}$ are given explicitly as functions of  $p, n,f$ and $\Omega$.  To establish the lower bound, we approximate the solution \( u \) to \eqref{HJ-S} by its sup-convolution \( u_\theta \), which is semiconvex. This enables a more direct estimate of the error \( u^\varepsilon - u_\theta \). However, the domain of  \( u_\theta \) must be dealt with carefully, as the inf-convolution alters the region where the equation remains valid. These steps are carried out in Subsection \ref{subsection:lowerbound}, where the lower bound is obtained in Proposition \ref{prop:LowerBound}.  For the upper bound, we establish a sharp estimate on $u^{\ve}$ in Proposition  \ref{lem:ZeroData}, refining the constant using a modified distance function. This improves upon existing bounds found in previous literature, such as \cite[Remark 4.5]{armstrong_viscosity_2015}.

It is to be noted that the rate $\mathcal{O}(\varepsilon^{1/2})$ in Theorem \ref{thm:NonnegativeZeroBoundaryIntro} naturally arises from the doubling of variables method, a standard technique in the study of viscosity solution (see \cite{bardi_optimal_1997, crandall_two_1984, fleming_convergence_1961}). In fact, there is an example in \cite{qian_optimal_2024} where this rate is optimal. 
However, when the data  \( f \) is a nonnegative semiconcave function with a compact support, the upper bound in Theorem \ref{thm:NonnegativeZeroBoundaryIntro} can be  improved. Indeed, set
\[
\alpha_p~\doteq~ \frac{p-2}{p-1}~\in~ (0,1)\qquad\forall p>2,
\]
we prove the following theorem.
\begin{theorem}\label{semi-im} For every $p>2$, assume that $\Omega\subset\R^n$ is an open and bounded subset of $\R^n$ with $\mathcal{C}^2$ boundary, $f\in \mathcal{C}_c(\overline{\Omega})$ is nonnegative and semiconcave in $\Omega$.

Then for every $\ve>0$ sufficiently small, it holds that
\bel{ke}
\max_{x\in\overline{\Omega}}\{u^{\ve}(x)-u(x)\}~\leq~\left({1\over\lambda}+\frac{2^{\alpha_p}}{\alpha_p}\right)\cdot \ve^{1-\frac{\alpha_p}{2}}~.
\eeq
\end{theorem}

For $p>2$, the solution $u$ may not be semiconcave for a semiconcave data $f$ since the  Legendre transform $L$  is no longer semiconcave; consequently,  the method in \cite{han_remarks_2022} does not apply. However, if $f\in \mathcal{C}_c(\overline{\Omega})$ is nonnegative and semiconcave with a semiconave constant ${\bf c}_f$ in $\Omega$, then by  Lemma \ref{semi-c}, the corresponding solution is still semiconcave with a semiconcavity constant ${\bf c}_{f}$. This holds because every optimal trajectory of (\ref{O-C}) for $x\in \Omega$ can stop upon reaching the support of $f$. As a result, estimate (\ref{ke}) follows from Proposition \ref{lem:ZeroData} via a standard argument. As a consequence, Corollary \ref{sf} shows that (\ref{ke}) can be obtained for a class of nonnegative data $f\in C^2(\Omega)$ satisfying  
\[
f(x)~=~0,\qquad Df(x)~=~0~~~\mathrm{on}~\partial\Omega~.
\]
In this case,  $f$ can be approximated by a sequence of semiconcave functions with compact support, and  Theorem \ref{semi-im} can be applied.  However, the upper estimates in Theorems \ref{thm:NonnegativeZeroBoundaryIntro} and \ref{semi-im} are yet to be proved for general Lipschitz data $f$ and remain open. When  the maximizer of $f$ lies in $\Omega$, the optimal trajectory may touch the boundary and then re-enter the domain. In such cases, a different approach is needed.  Some recent works related to convergence rates for Hamilton-Jacobi equations with vanishing viscosity can be found in \cite{chaintron_optimal_2025, cirant_convergence_2025,{qian_optimal_2024}, {tran_adjoint_2011}}, while \cite{han_quantitative_2025} addresses convergence rates in the problem with state-constraints.

The remainder of this paper is organized as follows. Section 2 is a review of key definitions, such as semiconcavity and viscosity solutions, along with known estimates for problems having state constraints. We also present a new, simpler proof of a local Lipschitz estimate for solutions to \eqref{HJ-ve}, leveraging the explicit form of the Hamiltonian. In Section 3, we establish upper and lower bounds for \( u^\varepsilon - u \). The lower bound uses sup-convolution and boundary properties of \( u^\varepsilon \), using a different approach from the \( p \leq 2 \) case in \cite{han_remarks_2022}. The upper bound adapts ideas from \cite{han_remarks_2022} while introducing a more delicate barrier construction suited to the less understood boundary behavior for the case of \( p > 2 \), yielding sharper estimates than those in \cite{armstrong_viscosity_2015}. Finally in subsection \ref{con-se}, we  improve the convergence rate under nonnegative semiconcave data with a compact support.

\section{Semiconcave functions  and  viscosity solutions}

\subsection{Semiconcave functions}
Given an open set $\Omega\subset\R^n$, we say that the function $f: \Omega\to \R$ is semiconcave with  a linear modulus in $\Omega$ and has a semiconcavity constant $K$ if for all $x,y\in \Omega$ such that the line segment $[x,y]=\{(1-t)x+ty:t\in [0,1]\} \subset \Omega$ and for all $\lambda \in [0,1]$, it holds that 
\begin{align}\label{eq:definitionSemiconcaveLmabda}
    \lambda f(x) + (1-\lambda)f(y) - f\big(\lambda x + (1-\lambda) y \big)~\leq~ \lambda(1-\lambda) \cdot \frac{K|x-y|^2}{2}~.
\end{align}

The following lemma is standard; for instance, see \cite{cannarsa_semiconcave_2004}.
\begin{lemma}
Assume that $f:\Omega\to\R$ is semiconcave  with a semiconcavity constant $K>0$. Then  $f$ is twice differentiable almost everywhere in $\Omega$ and 
\bel{vis-semi}
-D^2f ~\geq~- K\cdot\mathbf{I}_n
\eeq
in the viscosity sense.
\end{lemma}

The distance functions to $\partial\Omega$ and $\overline{\Omega}$ are defined by
\[
d_{\partial\Omega}(x)~=~\min_{y\in \partial\Omega} |y-x|,\qquad d_{\overline{\Omega}}(x)~=~\min_{y\in \overline{\Omega}} |y-x|\qquad\forall x\in \R^n.
\]
For every $\delta>0$, let
\bel{O-d}
\begin{split}
	\Omega_\delta ~\doteq~ \{x\in \Omega: d_{\partial\Omega}(x) > \delta\}, \qquad \text{and}\qquad 
	\Omega^\delta ~\doteq~ \{x\in \R^n: d_{\overline{\Omega}}(x) < \delta\}.
\end{split}
\eeq
It is known (see \cite[p. 669]{capuzzo-dolcetta_hamilton-jacobi_1990}, \cite{gilbarg_elliptic_2001}) that if  $\partial\Omega$ is $\mathcal{C}^2$ then both   \( \partial\Omega_\delta \) and \( \partial \Omega^\delta \) are also of class \( \mathcal{C}^2 \) for $\delta>0$ sufficiently small. Moreover, the distance function  $d_{\partial\Omega}$ is $\mathcal{C}^2$-smooth in a neighborhood of $\partial\Omega$. In this case, we shall denote by 
\begin{equation}\label{S-con}
    \delta_{\Omega}~\doteq~
    \frac{1}{3}
    \cdot \sup\left\lbrace 
    \delta\in 
    \left(0,\max_{x\in\Omega}d_{\partial\Omega}(x)\right): d_{\partial\Omega}~\text{is}~\mathcal{C}^2~\mathrm{smooth~in}~\Omega\backslash \Omega_{\delta}
    \right\rbrace 
    ~>~0.
\end{equation}
For every $x\in \Omega\backslash \Omega_{\delta_{4\Omega}}$, we have that 
\begin{align*}
    Dd_{\partial\Omega}(x)~=~\ds{x-\pi_{\partial\Omega}(x)\over |x-\pi_{\partial\Omega}(x)|}
\end{align*}
with $\pi_{\partial\Omega}(x)$ being the unique projection from $x$ to $\partial\Omega$. At each point $x_0\in \partial\Omega$, the inward normal vector to $\Omega$ at $x_0$ is given by $\mathbf{n}_{\partial\Omega}(x_0) = Dd_{\partial\Omega}(x_0)$. In addition, it holds that
\bel{in-n}
d_{\partial\Omega}(x_0 + s\mathbf{n}_{\partial\Omega}(x_0))~=~s\qquad\forall s\in [0,3\delta_\Omega],
\eeq
and 
\bel{bb}
\Omega\backslash\Omega_{\delta}~=~\bigcup_{x_0\in \partial\Omega}\left\{x_0+s\cdot\mathbf{n}_{\partial\Omega}(x_0):s\in [0,\delta], x_0\in \partial\Omega\right\},\qquad \delta\in [0,3\delta_\Omega].
\eeq

\begin{lemma} Assume that $\Omega\subset\R^n$ is open, bounded with a $\mathcal{C}^2$ boundary. Then  it holds that
\bel{De-d}
\Delta d_{\partial\Omega}(x)~\geq~- {n\over \delta_\Omega},\qquad x\in \Omega\backslash \Omega_{2\delta_{\Omega}}.
\eeq
\end{lemma}
{\bf Proof.} For every $x\in  \Omega\backslash \Omega_{2\delta_{\Omega}}$, we have 
\[
d_{\partial\Omega}(x)~=~2\delta_{\Omega}-d_{\partial\Omega_{2\delta_{\Omega}}}(x).
\]
By  the definition of $\delta_{\Omega}$ in (\ref{S-con}), the set $\Omega_{2\delta_{\Omega}}$ satisfies an interior sphere condition with radius $\delta_{\Omega}$. In particular, from \cite[Proposition 2.2.2]{cannarsa_semiconcave_2004} the function $x\mapsto d_{\partial \Omega_{2\delta_{\Omega}}}(x)$ is semiconcave with semiconcavity constant $1/\delta_\Omega$. As a consequence, for every $x\in  \Omega\backslash \Omega_{2\delta_{\Omega}}$, the matrix 
\begin{equation*}
	\ds D^2d_{\partial \Omega_{2\delta_{\Omega}}}(x)-{\mathrm{{\bf I}}_n\over \delta_{\Omega}}
\end{equation*}
is negative definite. Hence, we derive 
\[
\Delta d_{\partial\Omega_{2\delta_{\Omega}}}(x)~\leq~\mathrm{trace}\left({\mathrm{{\bf I}}_n\over \delta_{\Omega}}\right)~=~{n\over \delta_{\Omega}},
\]
which yields (\ref{De-d}).
\endproof

\subsection{Viscosity solutions to Hamilton-Jacobi equations with state constraints}
Given an open and bounded set $\Omega\subset\R^n$, let $H:\overline{\Omega}\times \R^n\to \R$ be a continuous Hamiltonian. For every $\ve\geq 0$,   we consider   the Hamilton-Jacobi equation
\begin{equation}\label{eq:pde}
	\lambda u(x) + H(x,Du(x)) - \varepsilon \Delta u(x)~=~ 0 \quad \text{in}~\Omega.
\end{equation}
\begin{definition} Let $\Omega\subset\R^n$ be an open subset. 
\begin{itemize}
\item[(i)] We say that an upper semicontinuous function $u:\Omega\to\R$ is a viscosity subsolution to \eqref{eq:pde} in $\Omega$, if for any $x_0\in \Omega$ and $\varphi\in C^2(\Omega)$ (or $\varphi\in C^1(\Omega)$ for $\ve=0$) such that $u-\varphi$ has a local maximum at $x_0$, it holds that 
\begin{equation*}
	\lambda u(x_0) + H(x_0,D\varphi(x_0)) - \varepsilon \Delta \varphi(x_0)~ \leq~ 0. 
\end{equation*}
\item[(ii)] We say that a lower semicontinuous function $v:\overline{\Omega}\to\R$ is a viscosity supersolution to \eqref{eq:pde} on $\overline{\Omega}$, if for any $x_0\in \overline{\Omega}$ and $\varphi\in C^2(\overline{\Omega})$  (or $\varphi\in C^1(\Omega)$ for $\ve=0$)  such that $u-\varphi$ has a local minimum at $x_0$, it holds that
\begin{equation*}
	\lambda v(x_0) + H(x_0,D\varphi(x_0)) - \varepsilon \Delta \varphi(x_0)~\geq~ 0. 
\end{equation*}
\end{itemize}
If $u:\overline{\Omega}\to\R$ is both a subsolution in $\Omega$, and a supersolution on $\overline{\Omega}$, we say that $u$ is a viscosity solution to a problem with state constraints, and satisfies 
\begin{equation}\label{eq:StateConstraintPDE}
\begin{cases}
	\begin{aligned}
		\lambda u(x) + H(x,Du(x)) - \varepsilon \Delta u(x) &~\leq~ 0  &&\text{in}~\;\Omega \\[2mm]
		\lambda u(x) + H(x,Du(x)) - \varepsilon \Delta u(x) &~\geq~ 0  &&\text{on}~\;\overline{\Omega}
	\end{aligned}
\end{cases}~.
\end{equation}
\end{definition}
Assume that  \( p \mapsto H(x,p) \) is convex with Legendre transform $L:\overline{\Omega}\times\R^n\to\R$ defined by \eqref{L}. Then the solution to the problem with state constraints can be expressed as the value function of an optimal control problem with running cost \( L(x,v) \), constrained to trajectories in \( \overline{\Omega} \) (see \cite{soner_optimal_1986}). The boundary condition with state constraints, which is an implicit condition, is often called the natural boundary condition in optimal control theory. Specially, in the case of  $\ve=0$, the unique viscosity solution of (\ref{eq:StateConstraintPDE}) is represented by \eqref{O-C}. 

\begin{remark} In general, if \( u \in C(\overline{\Omega}) \) is a solution to \eqref{eq:StateConstraintPDE}, its restriction \( v = u|_{\Omega'} \) for \( \Omega' \subset \Omega \) may not satisfy the condition of having state constraints on \( \Omega' \) (see \emph{\cite{kim_state-constraint_2020}} which studies error estimates for solutions on nested domains). 
\end{remark}

Next, we provide a self-contained proof of the comparison principle for \eqref{HJ-ve}.
We refer to \cite{armstrong_viscosity_2015, capuzzo-dolcetta_hamilton-jacobi_1990, crandall_users_1992, soner_optimal_1986} for related results. 

\begin{proposition}\label{prop:Comparison}
Assume that $\Omega\subset\R^n$ is an open and bounded subset of $\R^n$ with $\mathcal{C}^2$ boundary, and $H$ is of the form \eqref{H}, i.e.,
$$H(x,\xi) = |\xi|^p - f(x),\qquad (x,\xi)\in \overline{\Omega}\times \R^n.$$
Then, for every pair $(u,v)\in {\mathcal{C}}(\overline{\Omega})\times \mathcal{C}(\overline{\Omega})$, where $u$ is a viscosity subsolution in $\Omega$  of \eqref{eq:pde} with $f=f_u\in \mathcal{C}(\overline{\Omega})$ and $v$ is a supersolution on $\overline{\Omega}$ of \eqref{eq:pde}  with $f=f_v\in \mathcal{C}(\overline{\Omega})$,
it holds that 
\bel{cp}
	\sup_{\overline{\Omega}} (u-v)~\leq~\frac{1}{\lambda}
	\cdot 
	\sup_{\overline{\Omega}}\, (f_u-f_v)^+. 
\eeq
\end{proposition}
{\bf Proof.} Observe that  $\tilde{u} = u - \ds \frac{1}{\lambda}\cdot \sup_{\overline{\Omega}}\, (f_u-f_v)^+$ is a subsolution to 
\begin{align}\label{eq:svvv}
	\lambda w+ |Dw|^p - f_v - \varepsilon \Delta w~\leq~ 0 \qquad\text{in}\;\Omega.
\end{align}
Since $v$ is a supersolution of \eqref{eq:svvv} on $\overline{\Omega}$, for any $\psi\in \mathcal{C}^2(\overline{\Omega})$, at any local minimum $x_0\in \partial\Omega$ of $v-\psi$ in $\overline{\Omega}$
\begin{align*}
	& \max
		\left\lbrace 
		\lambda v(x_0) + |D\psi(x_0)|^p - f(x_0) - \varepsilon\Delta \psi(x_0)  	
		,
		v(x_0) - \psi(x_0)	
		\right\rbrace   \\
	&\qquad\qquad\qquad\qquad\qquad 	\geq~ \lambda v(x_0) + |D\psi(x_0)|^p - f(x_0) - \varepsilon\Delta \psi(x_0)~\geq~ 0.
\end{align*}
According to \cite{barles_generalized_2004}, in the generalized Dirichlet boundary problem
\[
\begin{cases}
\lambda w + |Dw|^p - f_v - \varepsilon \Delta w~\leq~ 0 \qquad\text{in}\;\Omega,\\[2mm]
w~=~u|_{\partial\Omega}\qquad\text{in}\;\partial\Omega,
\end{cases}
\]
$u$ is a subsolution and $v$ is  a a supersolution. Therefore, by the comparison principle  \cite[Theorem 2.3]{tabet_tchamba_large_2010} (which based on \cite[Theorem 3.1]{barles_generalized_2004} and generalization in \cite[Section 5]{barles_generalized_2004}), we obtain
\begin{align*}
	\tilde{u}~\leq~ v \qquad\text{on}\;\overline{\Omega},
\end{align*}
which implies the desired result. 
\endproof

\subsection{Properties of the solutions to problems with state constraints}
Next, we state and prove some useful Lipschitz and H\"older regularity properties of viscosity solutions to \eqref{HJ-S} and \eqref{HJ-ve} for the superquadratic case $p>2$. To do so, we  recall that  
\bel{al-p}
\alpha_p~\doteq~\frac{p-2}{p-1}~\in~(0,1). 
\eeq
\begin{proposition}\label{prop:Holder} Assume that  $f\in \mathcal{C}(\overline{\Omega})$. Then the following hold:
\begin{itemize}
\item [(i).] If $v\in \mathrm{USC}(\Omega)$ is a bounded viscosity subsolution of \eqref{HJ-ve} in $\Omega$, then $v$ is uniformly H\"older continuous with the H\"older exponent $\alpha_p$.
\item [(ii).] If $f\in \mathrm{Lip}(\overline{\Omega})$, there exists a (maximal) solution $u^\varepsilon\in \mathcal{C}^2(\Omega)\cap \mathcal{C}^{0,\alpha_p}(\overline{\Omega})$ of \eqref{HJ-ve} which satisfies 
\begin{equation}\label{eq:NegativeBoundary}
	u^{\ve}~\geq~{1\over \lambda}\cdot {\min_{z\in \overline{\Omega}}f(z)} \qquad\mathrm{and}\qquad \liminf_{\Omega \ni x\to x_0} \frac{u(x)-u(x_0)}{|x-x_0|^{\alpha_p}} ~<~ 0
\end{equation} 
for all $x_0\in\partial\Omega$.
\end{itemize}
\end{proposition}
{\bf Proof.} Refer to \cite{barles_short_2010} or \cite[Theorem 2.7]{capuzzo-dolcetta_holder_2010} for a proof of (i). For (ii), the existence of a \(\mathcal{C}^2\) maximal solution, its continuous extension to \(\overline{\Omega}\) and \eqref{eq:NegativeBoundary} follow from {\cite[Theorem I.2]{lasry_nonlinear_1989}}. As a maximal solution, it must coincide with the one constructed via Perron's method in \cite[Theorem 4.2]{armstrong_viscosity_2015}. Since $\ds{1\over \lambda}\cdot {\min_{z\in \overline{\Omega}}f(z)}$ is a subsolution, Perron's method ensures that \(\lambda u^\varepsilon \geq {\min_{z\in \overline{\Omega}}f(z)}\). By (i), it follows that \( u^\varepsilon \) can be extended to \( u^\varepsilon \in \mathcal{C}^{0,\alpha_p}(\overline{\Omega}) \).
\endproof
\medskip

The following lemma is standard and holds for general continuous, coercive Hamiltonians; see \cite[Theorem 2.1]{kim_state-constraint_2020} for reference. We include the proof here as it follows directly from the oscillation of \( f \in \mathcal{C}(\overline{\Omega}) \).
\begin{lemma}\label{lem:BoundLipU} Under the same assumptions as in Proposition \ref{prop:Holder}, the solution $u$ of \eqref{HJ-S} is Lipschitz with 
\bel{Lip-u}
\|u\|_{{{\bf Lip}}}~\leq~ \left(\mathrm{osc}_\Omega f\right)^{1/p}.
\eeq
\end{lemma}
{\bf Proof.} Observe that the constant function $v\equiv\ds\lambda^{-1}\cdot{\inf_{z\in \Omega} f(z)}$ is a viscosity subsolution to $\lambda u +|Du|^p - f(x)\leq 0$ in $\Omega$. By comparison principle, we have $\lambda u(x) \geq \ds\inf_{z\in \Omega} f(z)$ for $x\in \Omega$, which in turns implies that for every $\xi\in D^+u(x)$ and $x\in \Omega$,
\[
	\left(\inf_{z\in \Omega} f(z)\right) - f(x) + |\xi|^p~\leq~\lambda u(x) + |\xi|^p - f(x)~\leq~0 	.
\]
Hence, we have $|\xi|^p \leq f(x) - \inf_{z\in\Omega} f(z) \leq \mathrm{osc}_\Omega f$,
which yields that $u$ is Lipschitz on $\overline{\Omega}$ with Lipschitz constant $ \left(\mathrm{osc}_\Omega f\right)^{1/p}$.
\endproof

Next, we establish a (local) gradient bound for $u^\varepsilon$. Such an estimate is seen in \cite{armstrong_viscosity_2015} and \cite[Theorem IV.1]{lasry_nonlinear_1989}. 
The classical approach to such a local Lipschitz estimate, i.e., a bound on $Du^\varepsilon$, is the so-called Bernstein's method which requires differentiating the equation either in the classical sense, or within the viscosity sense via the doubling variable method; see \cite[Lemma 3.2]{armstrong_viscosity_2015} where an explicit constant is provided in terms of the sup norm of $u^\varepsilon$. 

We give an explicit Lipschitz constant depending on $\mathrm{osc}_\Omega f$ in the following lemma, and we do not use Bernstein method in the classical way, i.e., we do not differentiate the equation, and the constant we obtain is \emph{sharp} in the sense that it matches the asymptotic behavior of the gradient $Du^\varepsilon$ near the boundary for $1<p\leq 2$, as in \cite{alessio_asymptotic_2006}.

\begin{lemma}\label{lem:C0} Assume that $p>2$,  $f\in \mathcal{C}(\overline{\Omega})$. Let $u^\varepsilon\in \mathcal{C}^2(\Omega)$ be the maximal viscosity solution to \eqref{HJ-ve}. Then for every $x\in \Omega$, it holds that
\begin{equation}\label{sha}
	|Du^\varepsilon(x)|~\leq~  \left({p\over p-1}\cdot \mathrm{osc}_{\Omega}f+\left(\frac{\varepsilon}{d_{\partial\Omega}(x)}\right)^{p/(p-1)}\right)^{1/p}. 
\end{equation}
\end{lemma}
{\bf Proof.}  For clarity, we break the proof into two steps. 
\v

{\bf 1.} As in Lemma \ref{lem:BoundLipU},  since $v\equiv\ds\lambda^{-1}\cdot{\inf_{z\in \Omega} f(z)}$ is a viscosity subsolution to $\lambda u + |Du|^p - f(x) - \varepsilon \Delta u \leq 0$ in $\Omega$, it holds that
\bel{oc}
\lambda u^{\ve}(x)~\geq~\inf_{z\in \Omega}f(z)\quad\forall x\in \Omega.
\eeq
Fix  $x\in \Omega$. For every $r\in (0,d_{\partial\Omega}(x))$ such that $B_r(x)\subset\subset \Omega$, we consider the scaling function $v:B_1(0)\to\R$ defined by 
\begin{equation*}
	v(y)~=~ \frac{1}{r} \cdot u^\varepsilon\big( x + r y \big),\qquad  y\in B_1(0).
\end{equation*}
Note that, $v\in \mathcal{C}^2(\overline{B}_1(0))$ since $B_r(x)\subset\subset \Omega$. For every $y\in B_1(0),$ we compute 
\begin{equation*}
	|Dv(y)|~=~|Du^\varepsilon(x+ry)| \qquad\mathrm{and}\qquad \Delta v(y) = r \Delta u^\varepsilon(x + ry).
\end{equation*}
Hence, from (\ref{HJ-ve}) we deduce 
\begin{equation*}
	\lambda r v(y) +|Dv(y)|^p - {\varepsilon\over r}\cdot \Delta v(y)~=~ f(x+ry),
\end{equation*}
and then (\ref{oc}) yields 
\bel{stk}
|Dv(y)|^p -  {\varepsilon\over r}\cdot \Delta v(y)~\leq~\mathrm{osc}_{\Omega}f,\qquad y\in B(0,1).
\eeq

\n {\bf 2.} Next, for every $\delta>0$, let $\psi_{\delta}:  \overline{B_1(0)}\to [0,1]$ be smooth such that 
\[
|D\psi_{\delta}|~\leq~1+\delta,\qquad \psi_{\delta}(0)~=~1,\qquad \psi_\delta(x)~=~0\quad\forall  x\in\partial B_1(0).
\]
Since $v\in \mathcal{C}^2(\overline{B}_1(0))$, the map $x\mapsto \big|\psi_\delta(x) Dv(x)\big|$ achieves a maximizer $x_0\in B_1(0)$ over $\overline{B_1(0)}$. If $\big|\psi_\delta(x_0) Dv(x_0)|=0$ then 
\[
|Du^\varepsilon(x)|~=~|Dv(0)|~=~ |\psi_\delta(0)Dv(0)|~\leq~\big|\psi_{\delta}(x_0) Dv(x_0)\big|~=~0.
\]
Otherwise, we have $D\big(\psi_\delta(x_0)Dv(x_0)\big) = 0$, which yields
\[
D^2v(x_0)~=~- {1\over \psi_{\delta}(x_0)}\cdot D\psi_{\delta}(x_0)\otimes Dv(x_0)~,
\]
and this in particular yields 
\[
|\Delta v(x_0)|~\leq~ {|D\psi_{\delta}(x_0)|\over \psi_{\delta}(x_0)}\cdot |Dv(x_0)|~.
\]
Recalling (\ref{stk}), we get
\[
|Dv(x_0)|^p- {\ve\over r}\cdot {|D\psi_{\delta}(x_0)|\over \psi_{\delta}(x_0)}\cdot |Dv(x_0)|~\leq~\mathrm{osc}_{\Omega}f~.
\]
Since $p>2$, $0\leq \psi\leq 1$ and $|D\psi_\delta| \leq 1 + \delta$, we have
\[
|\psi_{\delta}(x_0)Dv(x_0)|^p-{(1+\delta)\ve\over r}\cdot |\psi_{\delta}(x_0)Dv(x_0)|~\leq~\mathrm{osc}_{\Omega}f~.
\]
Hence, by Young's inequality we derive 
\[
{p-1\over p}\cdot \left(|\psi_{\delta}(x_0)Dv(x_0)|^p- \left((1+\delta)\ve\over r\right)^{p/(p-1)}\right)~\leq~\mathrm{osc}_{\Omega}f,
\]
and this yields 
\[
\begin{split}
|Du^\varepsilon(x)| &~=~|\psi_{\delta}(0)Dv(0)|\\
&~\leq~|\psi_{\delta}(x_0)Dv(x_0)|
		~\leq~ \left({p\over p-1}\cdot \mathrm{osc}_{\Omega}f+\left((1+\delta)\ve\over r\right)^{p/(p-1)}\right)^{1/p},
\end{split}
\]
for every $r\in (0,d_{\partial\Omega}(x))$ and $\delta>0$. Finally, talking $\delta\to 0^+$ and $r\to d_{\partial\Omega}(x)^-$, we achieve (\ref{sha}).
\endproof

\begin{corollary}\label{Co} Under the assumptions of Lemma \ref{lem:C0}, if \( f \equiv \text{constant} \), then 
\[
 |Du^\varepsilon(x)| ~\leq~  \left(\ve\over d_{\partial\Omega}(x)\right)^{1/(p-1)},\qquad \forall x\in \Omega.
\]
\end{corollary}

\section{Rate of convergence of $u^{\ve}\to u$} \label{sec:RateLipschitzData}
In this section, we provide a proof of Theorem \ref{thm:NonnegativeZeroBoundaryIntro}, which is divided into the following two subsections: a lower estimate in Proposition \ref{prop:LowerBound}, and an upper estimate in Proposition \ref{prop:upper}. Before proceeding, let us recall the following constants, which depend on $f$, $p$ and $\Omega$:
\bel{alphaf}
\|f\|_{{\bf {Lip}}}~\doteq~\sup_{x\neq y}{|f(x)-f(y)|\over |x-y|},\qquad \mathrm{osc}_\Omega f~\doteq~\left(\sup_{z\in \Omega}f(z)\right)-\left(\inf_{z\in \Omega}f(z)\right)~.
\eeq 

\subsection{Lower bound of $u^{\ve}-u$} \label{subsection:lowerbound}
Observe that by  \eqref{eq:NegativeBoundary},  \( u^\varepsilon - u \) attains a minimum over \( \overline{\Omega} \) at \( x_0 \in \Omega\). Using \( u^\varepsilon \) as a smooth test function, we derive 
\begin{align*}
	\min_{\overline{\Omega}}\big(u^\varepsilon-u\big)~=~u^\varepsilon(x_0) - u(x_0)~\geq~ {\varepsilon\over \lambda}\cdot  \Delta u^\varepsilon(x_0).
\end{align*}
Hence, if  \( u \) is semiconvex with constant \( K \), then 
\[
\min_{\Omega}\big(u^\varepsilon-u\big)~\geq~ {\varepsilon\over \lambda}\cdot  \Delta u^\varepsilon(x_0) ~\geq~{nK\ve\over\lambda}.
\]
More generally, we establish the following lemma.
\begin{lemma}\label{bi-in} For a given $g\in \mathcal{C}(\overline{\Omega})$, let $v\in \mathbf{Lip}(\overline{\Omega})$ be a viscosity subsolution of 
\bel{HJ-g}
\lambda v+|Dv|^p-g~=~0\qquad\mathrm{in}~~\Omega.
\eeq
If $v$ is semiconvex in $\Omega$ with a semiconvexity constant $K$ then 
\bel{u-v}
\min_{\overline{\Omega} } \big( u^{\ve}-v \big) ~\geq~{1\over \lambda}\cdot \left(-Kn\varepsilon+\min_{\overline{\Omega}}\{f-g\}\right).
\eeq
\end{lemma}
{\bf Proof.} Let $x_0\in \overline{\Omega}$ be a minimizer of $u^{\ve}-v$ over $\overline{\Omega}$. By (\ref{eq:NegativeBoundary}) and the Lipschitz continuity of $v$, it follows that $x_0\in\Omega$. Indeed, assuming $x_0\in\partial\Omega$ leads to a contradiction
\begin{equation}\label{eq:boundary}
    0~>~\liminf_{x\to x_0}\left(\frac{u^\varepsilon(x)-u^\varepsilon(x_0)}{|x-x_0|^{\alpha_p}} \right)~\geq~  \liminf_{x\to x_0}\left(\frac{v(x)-v(x_0)}{|x-x_0|^{\alpha_p}} \right) ~=~  0.
\end{equation}
By the semiconvexity property of $v$, we have that $\Delta u^{\ve}(x_0) \geq -Kn$. Since $v$ is a subsolution of (\ref{HJ-g}) and $u^{\ve}\in \mathcal{C}^2(\Omega)$ is the solution of \eqref{eq:pde}, using $u^{\ve}$ as a test function for (\ref{HJ-g})  at $x_0$, we obtain
\begin{align*}
    & \lambda v(x_0)+|Du^{\ve}(x_0)|^p- g(x_0)~\leq~0 , \\[2mm]
    & \lambda u^\varepsilon(x_0) + |Du^\varepsilon(x_0)|^p - f(x_0) - \varepsilon \Delta u^\varepsilon(x_0)~=~ 0.
\end{align*}
Hence, we derive 
\[
u^\varepsilon(x_0)-v(x_0)~\geq~ \frac{1}{\lambda}\cdot\Big(  \varepsilon \Delta u^\varepsilon(x_0)+ f(x_0)-g(x_0)\Big)~\geq~\frac{1}{\lambda}\cdot
\left(-Kn\varepsilon +\inf_{\overline{\Omega}}\left ( f-g\right)\right),
\]
which completes the proof.
\endproof
\v

However,  since \( u \) is not generally semiconvex, our approach is to apply a sup-convolution to obtain a semiconvex approximation \( u_\delta \). This approximation is only a subsolution in a smaller domain \( \Omega(\delta) \), so we introduce an additional approximation step to extend the validity back to the full domain \( \Omega \). This is carried out in the following proposition, where the last step is the correction of domains between $\Omega(\delta)$ and $\Omega$.

\begin{proposition}\label{prop:LowerBound} Assume that  $f$ is in  ${\bf Lip}(\overline{\Omega})$. Then for every   $u^\varepsilon$ and $u$ which are the unique viscosity solutions of \eqref{HJ-ve} and \eqref{HJ-S} respectively, it holds that
\begin{equation}\label{eq:PropConclusionLowerBound}
    \min_{x\in \overline{\Omega}} \big( u^\varepsilon(x) - u(x)\big) 
    \geq 
    -\underline{\mathbf{\Lambda}}\cdot \sqrt{\varepsilon}, 
\end{equation}
where the constant $ \underline{\mathbf{\Lambda}}$ is computed as 
\begin{align}\label{eq:LowerBoundAA}
	 \underline{\mathbf{\Lambda}} 
    ~\doteq~
    \frac{1}{\lambda}\left[n + \left(\frac{1}{2}\right)^{p-1}\cdot \mathbf{k}^{p+1}+ p\cdot \mathbf{C}_\Omega\cdot \mathrm{sup}_\Omega f\cdot\mathbf{k} + \frac{3}{2}\mathbf{k}^2
		  \right],
\end{align}
with 
\begin{equation*}
\mathbf{C}_\Omega~=~  \frac{3}{\delta_\Omega} + \Vert D^2 d_{\partial\Omega}\Vert_{\mathcal{C}\left(\Omega\backslash \Omega_{\delta_\Omega}\right)}
, \qquad 
\mathbf{k}~=~ 2\left(\mathrm{osc}_\Omega f + \Vert f\Vert_{\mathbf{Lip}}\right)^{1/p}.
\end{equation*}
\end{proposition}

{\bf Proof.} In view of Lemma \ref{bi-in}, for every $\theta\in (0,1)$ we shall construct a subsolution $u_{\theta}$ of (\ref{HJ-g}) with $g=g_{\theta}\in \mathcal{C}(\overline{\Omega})$ such that 
\begin{equation*}
    \ds\min_{\overline{\Omega}}(f-g_{\theta})~\geq~\mathcal{O}(\theta).
\end{equation*}
{\bf 1.} Let $\bar{f}:\R^n\to\R $ be the Lipschitz extension of $f$ defined by
\begin{equation*}
    \bar{f}(x)~ =~ \min_{y\in \overline{\Omega}} \Big\lbrace f(y) + \Vert f\Vert_{\mathbf{Lip}}\cdot|x-y| \Big\rbrace\quad\forall~ x\in \R^n,
\end{equation*}
such that 
\bel{ex-f}
\|\bar{f}\|_{\mathbf{Lip}}~=~\|f\|_{\mathbf{Lip}}\qquad\mathrm{and}\qquad \bar{f}~=~f~~\mathrm{in}~~\overline{\Omega},
\eeq
and 
\bel{ex-f1}
    \Vert f\Vert_{\mathrm{Lip}}\cdot d_{\overline{\Omega}}(x)+ \min_{\overline{\Omega}} f  ~\leq~\bar{f}(x)~\leq~ \max_{ \overline{\Omega}} f + \Vert f\Vert_{\mathrm{Lip}}\cdot  d_{\overline{\Omega}}(x), \qquad x\in \R^n.
\eeq
For any \( \theta > 0 \), consider the following open set in \( \mathbb{R}^n \) containing \( \Omega \):
\begin{align}\label{Omega-t}
	\Omega^{\mathbf{k}_f\theta}
	~\doteq~ 
	\left\{x\in \R^d: d_{\overline{\Omega}}(x)<\mathbf{k}\theta\right\}
	\qquad\mathrm{with}\qquad {\bf k} 
	~\doteq~
	2\left(\mathrm{osc}_\Omega f + \Vert f\Vert_{\mathbf{Lip}}\right)^{1/p}~.
\end{align}
We denote by  $v_{\theta}$  the viscosity solution of 
\bel{HJ-theta}
\begin{cases}
\lambda v(x)+ |D v(x)|^{p}~=~\bar{f}(x),\qquad x\in \Omega^{\mathbf{k}\theta}\\[2mm]
\lambda v(x)+ |D v(x)|^{p} ~\geq~ \bar{f}(x), \qquad x\in \partial\Omega^{\mathbf{k}\theta}
\end{cases}~.
\eeq
From  the representation formula (\ref{O-C}) and  (\ref{ex-f}), it holds that
\bel{nnq}
v_{\theta}(x)~\leq~u(x)~\leq~{1\over \lambda}\cdot \sup_{\overline{\Omega}} f\qquad\forall x\in \overline{\Omega}.
\eeq
By Proposition \ref{prop:Comparison}, (\ref{ex-f1}) and (\ref{Omega-t}), we have 
\begin{align}
    \|v_{\theta}\|_{{\mathrm{Lip}}}
    &~\leq~
    \Big(\mathrm{osc}_{\Omega^{\mathbf{k}\theta}} \bar{f}\Big)^{1/p} 
    ~\leq~ 
    \left(\mathrm{osc}_\Omega f+\|f\|_{{\mathbf{Lip}}}\cdot\mathbf{k}\theta\right)^{1/p} 
   ~ <~
   \left(\mathrm{osc}_\Omega f + \Vert f\Vert_{\mathbf{Lip}}\right)^{1/p}~ =~ \frac{\mathbf{k}}{2}
    \label{eq:LipVtheta}
\end{align}
provided \( \mathbf{k}\theta < 1 \).
\medskip

\n {\bf 2.} 
Next, we define the sup-convolution $u_\theta:\R^d \to \R$ of $v_{\theta}$ by
\begin{equation}\label{eq:supconvolutionu}
   u_\theta(x)~\doteq~ \sup_{y\in \overline{\Omega}} \left\lbrace v_{\theta}(y) - \frac{|x-y|^2}{2\theta} \right\rbrace, \qquad x\in \R^d. 
\end{equation}
The map $x\mapsto u_{\theta}(x)$ is semiconvex with a semiconvexity constant $1/\theta$, i.e., $-D^2u_\theta \geq (1/\theta)\;\mathbf{I}_n$ in the viscosity sense, and satisfies
\begin{equation}\label{eq:compareUthetaVtheta}
	u_{\theta}(x) - \mathbf{k}^2\theta 
	= 	
	u_{\theta}(x) - 4\Vert v_\theta\Vert_{\mathbf{Lip}}^2\cdot\theta 
	~\leq~
	v_{\theta}(x)~\leq~u_\theta(x),\quad\forall x\in \overline{\Omega}^{\mathbf{k}\theta}.
\end{equation}
In addition, $u_{\theta}$ is a viscosity subsolution of 
\[
\lambda u_\theta(x)+ |D u_\theta(x)|^{p}-g_{\theta}(x)~=~0,\qquad x\in U^{\mathbf{k}\theta},
\]
where $g$ and $U^{\mathbf{k}\theta}$ are defined by 
\begin{align*}
	g_{\theta}(x) &~=~ 
	f(x) 	+
	2\Vert \overline{f}\Vert_{\mathbf{Lip}(\Omega^{\mathbf{k}\theta})}
	\cdot \Vert v_\theta\Vert_{{\bf Lip}(\Omega^{\mathbf{k}\theta})}\cdot \theta, \\[3mm]
	U^{\mathbf{k}\theta}
	&~=~ \left\{x\in \Omega^{\mathbf{k}\theta}: \ds\mathrm{argmax}_{y\in \overline{\Omega}^{\mathbf{k}\theta}} \left\lbrace v_{\theta}(y) - \frac{|x-y|^2}{2\theta}\right\rbrace  \cap \Omega^{\mathbf{k}\theta} \neq \emptyset \right\} .
\end{align*}
We observe from \eqref{eq:LipVtheta} that 
\begin{align}\label{eq:fg}
	\min_{\overline{\Omega}}(f-g_\theta)~ \geq ~-\Vert f\Vert_{\mathrm{Lip}}\cdot
	\mathbf{k}\cdot \theta. 
\end{align}
For every $\theta>0$  such that  ${\bf k}\theta < 1$, we claim that $\Omega\subset U^{\mathbf{k}\theta}$. Indeed, for every $x\in \Omega$ and  
\begin{align*}
	y_x~ \in~ \mathrm{argmax}_{y\in \overline{\Omega}^{\mathbf{k}\theta}} 
	\left\lbrace v_{\theta}(y) - \frac{|x-y_x|^2}{2\theta}\right\rbrace.
\end{align*}
Then from (\ref{eq:compareUthetaVtheta}) we have 
\[
\frac{|x-y_x|^2}{2\theta}~=~ v_\theta(y_x) - u_\theta(x)  ~\leq~ v_\theta(y_x) - v_\theta(x) ~\leq~ \Vert v_\theta\Vert_{\mathbf{Lip}}\cdot |x-y_x|.
\]
In particular, we deduce from (\ref{eq:LipVtheta}) that 
\[
|x-y_x|~\leq~ 2\Vert v_\theta\Vert_{\mathbf{Lip}}\cdot \theta ~<~ \mathbf{k}\theta. 
\]
By the definition of $\Omega^{\mathbf{k}\theta}$ and $U^{\mathbf{k}\theta}$, for $x\in \Omega$ we get $y_x\in \Omega^{\mathbf{k}\theta}$, which  implies that  $u_\theta$ is a viscosity subsolution to
\begin{equation*}
	\lambda u_\theta + |D u_\theta|^p - g_\theta(x)~ =~ 0\qquad\text{in}~~\Omega. 
\end{equation*}
Hence, by Lemma \ref{bi-in} and \eqref{eq:fg} we obtain
\begin{align*}
	&\min_{\overline{\Omega}} (u^\varepsilon - u_\theta) 
	~\geq~ 
	-\frac{1}{\lambda}\cdot
	\left(
		\frac{n\varepsilon}{\theta} + \Vert f\Vert_{\mathbf{Lip}}\cdot \mathbf{k}\cdot\theta 
	\right)
\end{align*}.
As a consequence, from \eqref{eq:compareUthetaVtheta} we have
\begin{align}
    \min_{\overline{\Omega}} ( u^\varepsilon - v_\theta) 
    &~\geq~
    -\frac{1}{\lambda}\cdot
	\left(
		\frac{n\varepsilon}{\theta} 
		+ 
		\Vert f\Vert_{\mathbf{Lip}}\cdot \mathbf{k}\cdot\theta 
		+ \mathbf{k}^2\theta		
	\right)
\label{eq:uepsvtheta}.
\end{align}

{\bf 3.} To complete the proof, we  establish a lower  bound on $\ds\min_{\overline{\Omega}} \big(v_\theta - u\big)$. For $\theta > 0$ sufficiently small such that 
\begin{equation}\label{eq:conditionTheta}
	\delta~\doteq  ~\mathbf{k}\theta ~<~  \min \left\lbrace 1, \delta_\Omega/3 \right\rbrace,
\end{equation}
let $\eta_{\delta}:\R\to\R$ be a decreasing $\mathcal{C}^2$-smooth function such that 
\bel{eta-d}
\frac{3}{\delta_\Omega}~\leq~ \eta_\delta'~<~0,\qquad\quad \eta_{\delta}(s)~=~\begin{cases}
1,&~~ s\in [0,\delta]\\[2mm]
0,&~~ s\geq \delta_{\Omega}
\end{cases}~.
\eeq
Recalling that $\Omega_{\delta}=\{x\in \Omega:d_{\partial\Omega}(x)<\delta\}$, we consider $T_{\delta}:\overline{\Omega}\to \overline{\Omega}^\delta$ defined by
\begin{align*}
	T_{\delta}(x)~=~x - \delta \eta_{\delta}\left(d_{\partial\Omega}(x)\right)\cdot D d_{\partial\Omega}(x), \qquad x\in \overline{\Omega}. 
\end{align*}
Recalling (\ref{S-con}), we  have
\[
D d_{\partial\Omega}(x)~=~{x-\pi_{\partial{\Omega}}(x)\over |x-\pi_{\partial{\Omega}}(x)|}\quad\forall x\in \Omega\backslash \Omega_{\delta_{\Omega}}
\]
and this yields 
\[
 T_{\delta}(x)~=~\begin{cases}x,&\qquad  x\in \overline{\Omega}_{\delta_\Omega}\\[2mm]
\ds x-\delta\cdot \eta_{\delta}(d_{\partial \Omega}(x)) \cdot {x-\pi_{\partial{\Omega}}(x)\over |x-\pi_{\partial{\Omega}}(x)|},&\qquad x\in \Omega_{\delta}\backslash\Omega_{\delta_{\Omega}} \\[2mm]
\ds x-\delta \cdot {x-\pi_{\partial {\Omega}}(x)\over |x-\pi_{\partial {\Omega}}(x)|},&\qquad x\in \overline{\Omega} \backslash \Omega_\delta
\end{cases}~~~~.
\]
In particular,  the map $T_{\delta}:\overline{\Omega}\to \overline{\Omega}^\delta$ is bijective and satisfies 
\bel{eT}
    T_{\delta}(\partial\Omega)
    ~=~\partial\Omega^{\delta},
        \quad 
    d_{\overline{\Omega}}(T_{\delta}(x))~\leq~\big|T_{\delta}(x) - x\big|~\leq~\delta,
        \quad 
    \sup_{x\in \overline{\Omega}}\{|D T_{\delta}(x) - \mathbf{Id}|\} ~\leq~\mathbf{C}_\Omega\delta
\eeq
where  the  constant ${\bf C}_{\Omega}>0$ is computed by 
$$
\mathbf{C}_\Omega~\doteq~\ds  \frac{3}{\delta_\Omega}+\sup_{x\in\Omega\backslash\Omega_{\delta_{\Omega}}} \big|D^2 d_{\partial\Omega}(x)\big|~.
$$
Consider  $\tilde{v}_{\delta}:\overline{\Omega}\to\R$ such that 
\begin{equation*}
	\tilde{v}_{\delta}(x)~=~ v_{\theta}\big(T_{\delta}(x)\big)\qquad\forall x\in \overline{\Omega}. 
\end{equation*}
By (\ref{eq:LipVtheta}), (\ref{eT}), and as $\mathbf{k}\theta < 1$, we have 
\begin{align*}
    |\tilde{v}_{\delta}(x)-v_{\theta}(x)|
    ~\leq~ 
    \left|v_{\theta}(T_\delta(x))-v_{\theta}(x)\right|
    ~\leq~ \frac{\mathbf{k}\delta}{2}~, 
\end{align*}
which yields 
\bel{v-u}
    \min_{\overline{\Omega}}\big(v_{\theta}-u\big) ~\geq~
    \min_{\overline{\Omega}}\big( \tilde{v}_{\delta}-u\big) - \frac{\mathbf{k}\delta}{2}. 
\eeq
On the other hand, from \eqref{eT} we have 
\begin{align*}
	\sup_{x\in \overline{\Omega}}\big|\left(DT_\delta(x)\right)^{-1}\big|~\leq~ 1 + \frac{\mathbf{C}_\Omega \delta}{1-\mathbf{C}_\Omega \delta}~ =~ \frac{1}{1-\mathbf{C}_\Omega\delta}~~.
\end{align*}
Since $v_{\theta}$ is a viscosity solution of (\ref{HJ-theta}), $\tilde{v}_{\delta}$ is a viscosity supersolution of 
\begin{align*}
	\lambda \tilde{v}_\delta(x) + \big|(DT_\delta(x))^{-1}\cdot
	 D\tilde{v}_\delta(x)\big|^p ~=~ \bar{f}\big(T_\delta(x)\big), \qquad x\in \overline{\Omega}
\end{align*}
and thus, $\tilde{v}_\delta$ is a viscosity supersolution to 
\begin{align*}
	\lambda \tilde{v}_\delta(x) +\left(\frac{1}{1-\mathbf{C}_\Omega}\right)^p |D\tilde{v}_\delta(x)|^ p ~=~ \bar{f}\big(T_\delta(x)\big), \qquad x\in \overline{\Omega}. 
\end{align*}
In other words, $\tilde{v}_\delta$ is a viscosity supersolution to 
\begin{align*}
	\left(1-\mathbf{C}_\Omega \delta \right)^p \lambda \tilde{v}_\delta(x) + |D\tilde{v}_\delta(x)|^ p~ =~ \left(1-\mathbf{C}_\Omega \delta \right)^p \bar{f}\big(T_\delta(x)\big), \qquad x\in \overline{\Omega}. 
\end{align*}
Since $\left(1-\mathbf{C}_\Omega \delta \right)^p\leq 1$, 
$\tilde{v}_\delta$ is a viscosity supersolution to
\begin{align*}
	\lambda \tilde{v}_\delta(x) + |D\tilde{v}_\delta(x)|^ p~=~ \left(1-\mathbf{C}_\Omega \delta  \right)^p \bar{f}\big(T_\delta(x)\big), \qquad x\in \Omega. 
\end{align*}
Applying the standard comparison principle with the supersolution $\tilde{v}_\delta$ on $\overline{\Omega}$ and the subsolution $u$ in $\Omega$ (see \cite[Theorem III.1]{capuzzo-dolcetta_hamilton-jacobi_1990} for instance), we deduce that 
\begin{align*}
	\lambda\cdot \max_{\overline{\Omega}} (u - \tilde{v}_\delta) ^+ 
	&~\leq~ 
	\max_{x\in \overline{\Omega}} 
	\Big( 
	\left(1-\mathbf{C}_\Omega \delta  \right)^p 
	\bar{f}\big(T_\delta(x)\big) - f(x)
	\Big)^+\\
	&~=~ 
	\left(1-\mathbf{C}_\Omega \delta  \right)^p
	\cdot \max_{x\in \overline{\Omega}} 
	\left| \bar{f}(T_\delta(x)) - f(x) \right|  + \Big(1-\left(1-\mathbf{C}_\Omega \delta  \right)^p\Big)\cdot \sup_{\overline{\Omega}} f \\
	&~\leq~ (1-\mathbf{C}_\Omega\delta)^p \cdot \Vert f\Vert_{\mathbf{Lip}}\cdot \sup_{x\in \overline{\Omega}}|T_\delta(x) - x| 
	+ 
	p\cdot \mathbf{C}_\Omega \delta \cdot \sup_{\overline{\Omega}}f \\
	&~\leq~ \Big(\Vert f\Vert_{\mathbf{Lip}} + p\cdot \mathbf{C}_\Omega\cdot \sup_{\overline{\Omega}} f \Big)\cdot \delta~~.
\end{align*}
Using \eqref{ex-f}, \eqref{eT}, and Bernoulli's inequality: $(1-\kappa)^p\geq 1-p\kappa$ for all $p\geq 1$ and $\kappa \in [0,1]$, we conclude from \eqref{v-u} that 
\begin{align}\label{eq:vthetau}
	\min_{\overline{\Omega}} 
	\big( 
		 v_\theta - u
	\big) 
	~\geq~ 
        - \left[ 
        \frac{1}{\lambda}
        \Big(
            \Vert f\Vert_{\mathbf{Lip}} + p\cdot \mathbf{C}_\Omega\cdot \sup_{\overline{\Omega}} f 
        \Big) + \frac{\mathbf{k}}{2}
        \right]
    \cdot \delta ~~.
\end{align}
Combining \eqref{eq:uepsvtheta} and \eqref{eq:vthetau} with the facts \( \delta = \mathbf{k} \theta \) and $\lambda \in (0,1)$, we arrive at 
\begin{align}
	 \min_{\overline{\Omega}} \big( u^\varepsilon-u\big) 
	 &~\geq~ -\frac{1}{\lambda}
	 \left[
	 	\frac{n\varepsilon}{\theta} 
	 	+
		\left( 	
			2\Vert f\Vert_{\mathbf{Lip}}\cdot\mathbf{k} 
		 	+ 
	 		p \cdot \mathbf{C}_\Omega\cdot \sup_\Omega f\cdot \mathbf{k}
	 		+ 
	 		\frac{3}{2}\mathbf{k}^2
	 	\right)
	  \cdot \theta 
	 \right]  \label{eq:LowerInA}\\
	  &~\geq~ -\frac{1}{\lambda}
	 \left[
	 	\frac{n\varepsilon}{\theta} 
	 	+
		\left( 	
			\left(\frac{1}{2}\right)^{p-1} \cdot\mathbf{k}^{p+1}
		 	+ 
	 		p \cdot \mathbf{C}_\Omega\cdot \sup_\Omega f\cdot \mathbf{k}
	 		+ 
	 		\frac{3}{2}\mathbf{k}^2
	 	\right)
	  \cdot \theta 
	 \right].  \label{eq:LowerInB}
\end{align}
Finally, choosing  $\theta = \varepsilon^{1/2}$ in \eqref{eq:LowerInB}, we obtain \eqref{eq:PropConclusionLowerBound} with the constant given in \eqref{eq:LowerBoundAA}.
\endproof

\subsection{Upper bound of $u^{\ve}-u$ for Lipschitz data}\label{subsection:Upper}
In this subsection, we shall provide the upper estimate for $u^{\ve}-u$ in Theorem \ref{thm:NonnegativeZeroBoundaryIntro}. To this end, the next lemma establishes an upper bound for $u^{\ve}$ in the case of a constant function $f$,  offering an improvement over existing results in previous literature. For instance, in \cite[Lemma 4.3]{armstrong_viscosity_2015}, the author considered the case \( 1 < p \leq 2 \) and obtained a local bound of order \( \varepsilon^{\frac{1}{p-1}} = \varepsilon^{1-\alpha_p} \). Later, for \( p > 2 \), Remark 4.5 in the proof of \cite[Theorem 4.2]{armstrong_viscosity_2015} suggests a similar bound. In comparison, our estimate yields \( \varepsilon^{1-\frac{1}{2} \alpha_p} \), providing a refined bound in this setting.
 
\begin{lemma}\label{lem:ZeroData} For $p>2$, assume that  $f\equiv C_f$ is a constant function . Then for $\ve>0$ sufficiently small,  the unique viscosity solution $u^\varepsilon\in \mathcal{C}(\overline{\Omega})\cap \mathcal{C}^2(\Omega)$ of the problem \eqref{HJ-ve} with state constraints
satisfies
\begin{equation}\label{eq:EstimateUepsZerof}
	0~\leq~u^{\ve}(x)-{C_f\over \lambda}~\leq~\left(\ds{1\over\lambda}+\frac{2^{\alpha_p}}{\alpha_p}\right)\cdot \ve^{1-\frac{1}{2}\alpha_p}\quad\forall x\in  \overline{\Omega}.
\end{equation}
\end{lemma}

{\bf Proof.} Assume that  $f\equiv C_f$ is a constant function. Then the function $\tilde{u}^\varepsilon = u^\varepsilon - \frac{C_f}{\lambda}$ solves \eqref{HJ-ve} with zero data. Hence, we only need to establish (\ref{eq:EstimateUepsZerof}) for $f=0$. 
\medskip

\n {\bf 1.} It is clear that when $f=0$, the constant function $u\equiv 0$ is a classical  solution to  \eqref{HJ-S} and is also a classical subsolution to \eqref{HJ-ve}. Thus, by the standard comparison principle and Corollary \ref{Co}, we obtain that 
\bel{du-ve}
u^{\ve}(x)~\geq~0,\qquad |D u^\varepsilon(x)|~\leq~ \frac{\varepsilon^{1-\alpha_{p}}}{d_{\partial\Omega}(x)^{1-\alpha_p}},\qquad x\in \Omega.
\eeq
For every $\delta\in \left(0,\frac{1}{2}\delta_{\Omega}\right]$, let  $\eta_{\delta}:[0,\infty)\to [0,2\delta]$  be  a $\mathcal{C}^2$ increasing function   such that
\bel{eta}
\eta_{\delta}(s)~=~\begin{cases}
s,&\qquad s\in [0,\delta]\\[1mm]
2\delta,&\qquad s\in [2\delta_{\Omega},\infty)
\end{cases}
\eeq
and satisfies 
\bel{d-e}
0~\leq~\eta_{\delta}'(s)~\leq~1,\qquad  -{8\over \delta_{\Omega}}~\leq~ \eta_{\delta}''(s)~\leq~0,\qquad s\in [0,\infty).
\eeq
Introducing a $\mathcal{C}^2$-function  ${\bf d}_{\delta}:\overline{\Omega}\to [0,\infty)$ such that 
\bel{d-t}
{\bf d}_{\delta}(x)~=~\eta_{\delta} (d_{\partial\Omega}(x))\quad\forall x\in \overline{\Omega},
\eeq
we consider a  modification of  $u^{\ve}$ defined by 
\bel{w-ve}
w_{\delta}^\varepsilon(x)~=~ u^\varepsilon(x)+\ds \frac{1}{\alpha_p}\cdot\varepsilon^{1-\alpha_p}\cdot {\bf d}_{\delta}(x)^{\alpha_p}\qquad\forall x\in \overline{\Omega}.
\eeq

{\bf 2.} We first claim that $w^{\ve}_{\delta}$ admits a  global maximizer in $\Omega_{\delta}$. From (\ref{eta}) and (\ref{d-e}), we get
\[
w_{\delta}^\varepsilon(x)~=~ u^{\ve}(x)+
\ds\frac{1}{\alpha_p}\cdot\varepsilon^{1-\alpha_p}d_{\partial\Omega}(x)^{\alpha_p},\qquad x\in \Omega\backslash\Omega_{\delta}.
\]
For every   $x_0\in \partial\Omega$, let ${\bf n}_{\partial\Omega}(x_0)$ be the inner normal of $\Omega$ at $x_0$. By (\ref{in-n}),(\ref{eta}) and (\ref{d-e}), for every  $s\in [0,\delta]$, we have $x_0+s\cdot \mathbf{n}_{\partial\Omega}(x_0)\in \Omega\backslash\Omega_{\delta} $ and 
\[
w^\varepsilon_\delta(x_0+s\cdot \mathbf{n}_{\partial\Omega}(x_0))~=~\ds u^\varepsilon(x_0+s\cdot \mathbf{n}_{\partial\Omega}(x_0))+\ds\frac{1}{\alpha_p}\cdot\varepsilon^{1-\alpha_p} s^{\alpha_p}.
\]
Hence, using (\ref{du-ve}), we compute
\[
\begin{split}
{d\over ds}w^{\ve}_{\delta}(x_0+s\cdot \mathbf{n}_{\partial\Omega}(x_0))&~=~Du^{\ve}(x_0+s\cdot \mathbf{n}_{\partial\Omega}(x_0))\cdot {\bf n}_{\partial\Omega}(x_0)+{\varepsilon^{1-\alpha_p}\over s^{1-\alpha_p}}\\
&~\geq~ -{\varepsilon^{1-\alpha_p} \over  d_{\partial\Omega}^{1-\alpha_p}(x_0+s\cdot \mathbf{n}_{\partial\Omega}(x_0))}+{\varepsilon^{1-\alpha_p}\over s^{1-\alpha_p}}~=~0,
\end{split}
\]
which implies that  the map $s\mapsto w^\varepsilon_\delta(x_0+s\cdot \mathbf{n}_{\partial\Omega}(x_0))$ is non-decreasing in $[0,\delta]$, and  (\ref{bb}) implies that $w^{\ve}_{\delta}$ admits a  global maximizer in $\Omega_{\delta}$. \medskip

{\bf 3.} Next, picking any global maximizer $x_{\ve}\in \Omega_{\delta}$ of $w^{\ve}_{\delta}$ over $\overline{\Omega}$, we have 
\bel{dc}
D w^{\ve}_{\delta}(x_\varepsilon)~=~0\qquad\mathrm{and}\qquad \Delta w^{\ve}_{\delta}(x_{\ve})~\leq~0~.
\eeq
From \eqref{w-ve} and (\ref{dc}), we deduce
\[
|D u_{\ve}(x_\ve)|~=~ \left({\ve\over {\bf d}_{\delta}(x_{\ve})}\right)^{1-\alpha_p}\cdot |D{\bf d}_{\delta}(x_\ve)|,
\]
\[
\begin{split}
\Delta u^{\ve}(x_{\ve})&~\leq~-\ds \varepsilon^{1-\alpha_p}\cdot \Delta {\bf d}_{\delta}^{\alpha_p}(x_{\ve})~=~{\varepsilon^{1-\alpha_p}\over {\bf d}_{\delta}^{2-\alpha_p}(x_{\ve}) }\cdot \left[(1-\alpha_p) |D{\bf d}_{\delta}(x_\ve)|^2-  {\bf d}_{\delta}(x_{\ve})\Delta {\bf d}_{\delta}(x_{\ve})\right].
\end{split}
\]
Two cases are considered: 
\begin{itemize}
\item If $x_{\ve}\in \Omega_{2\delta_{\Omega}}$, then ${\bf d}_{\delta}(x_{\ve})=2\delta$, $D{\bf d}_\delta(x_\varepsilon) = 0$, $\Delta {\bf d}_{\delta}(x_\varepsilon) = 0$, and thus $\Delta u^\varepsilon(x_\varepsilon) = \Delta w^\varepsilon_\delta(x_\varepsilon)$ from \eqref{w-ve}. Therefore, together with \eqref{dc} we have
\[
u^{\ve}(x_{\ve})~=~{1\over \lambda}\cdot \left(\ve \Delta u^{\ve}(x_{\ve})-|D u^{\ve}(x_{\ve})|^p\right)~\leq~{\ve\Delta u^{\ve}(x_{\ve})\over\lambda}~=~{\ve\Delta w_{\delta}^{\ve}(x_{\ve})\over\lambda}~\leq~0.
\]
Hence, (\ref{du-ve}) yields $u^{\ve}(x_{\ve})=0$.
\item Otherwise if $x_{\ve}\in \Omega_{\delta}\backslash\Omega_{2\delta_\Omega}$, then $\delta \leq d_{\partial\Omega}(x_\varepsilon) \leq 2\delta_\Omega$. From \eqref{eta} and \eqref{d-t}, as $\eta$ is increasing on $[\delta, 2\delta_\Omega]$, we have
\[
	\delta~\leq~ {\bf d}_\delta(x_\varepsilon)~\leq~ 2\delta. 
\]
On the other hand, from the PDE for $u^\varepsilon$ we have
\bel{e1-uu}
\begin{split}
u^{\ve}(x_{\ve})&~=~{1\over \lambda}\cdot \left(\ve\cdot \Delta u^{\ve}(x_{\ve})-|D u^{\ve}(x_{\ve})|^p\right)\\
&~\leq~{\ve^{2-\alpha_p}\over \lambda{\bf d}_{\delta}^{2-\alpha_p}(x_{\ve})}\cdot \left[(1-\alpha_p)|D{\bf d}_{\delta}(x_\ve)|^2- {\bf d}(x_{\ve})\Delta {\bf d}_{\delta}(x_{\ve})\right].
\end{split}
\eeq
Recalling (\ref{eta}), we compute
\begin{align}
    |D{\bf d}_{\delta}(x_\ve)| 
    &~=~|\eta'_{\delta}(d_{\partial\Omega}(x_{\ve}))|~\leq~1,     
    \label{eq:estDBFd}
    \\
    \Delta {\bf d}_{\delta}(x_{\ve})
    &~=~\eta'_{\delta}(d_{\partial\Omega}(x_{\ve}))\cdot\Delta d_{\partial\Omega}(x_{\ve})+\eta''_{\delta}(d_{\partial\Omega}(x_{\ve}))~\geq~-{8+n\over\delta_\Omega},
    \label{eq:estDeltaBFd}
\end{align}
and (\ref{e1-uu}) yields 
\[
u^{\ve}(x_{\ve})~\leq~{\ve^{2-\alpha_p}\over \lambda\delta^{2-\alpha_p}}\cdot \left[1-\alpha_p+{(16+2n)\delta\over \delta_{\Omega}}\right]~\leq~{\ve^{2-\alpha_p}\over \lambda\delta^{2-\alpha_p}}
\]
by choosing $\delta$ small enough.
\end{itemize}
We observe that in both cases, we have ${\bf d}_\delta(x_\varepsilon) \leq 2\delta$.
Hence,  for every $x\in \Omega$, we have 
\[
\begin{split}
	u^\varepsilon(x) 
	&~\leq~ u^\varepsilon(x_\varepsilon) +\frac{1}{\alpha_p}\cdot \varepsilon^{1-\alpha_p}\cdot{\bf d}^{\alpha_p}_{\delta}(x_{\ve})~\leq~{\ve^{2-\alpha_p}\over \lambda\delta^{2-\alpha_p}}+ 2^{\alpha_p}\frac{1}{\alpha_p}\cdot \delta^{\alpha_p}\cdot \varepsilon^{1-\alpha_p}\\
	&~=~{1\over\lambda}\cdot\left({ \ve\over   \delta^{2-\alpha_p}}+{2^{\alpha_p}\lambda\over \alpha_p}\cdot \delta^{\alpha_p}\right)\cdot \ve^{1-\alpha_p}~~.\end{split}
\]
Finally, choosing $\delta = \sqrt{\varepsilon}$ we obtain (\ref{eq:EstimateUepsZerof}).
\endproof

\begin{proposition}\label{prop:upper}
For every $p>2$, assume that  $f\in {\bf Lip}(\overline{\Omega})$ is nonnegative and vanishes  on $\partial\Omega$.  Then 
\begin{equation}\label{eq:THMEstimateUepsConstantfUpper}
    \max_{x\in \overline{\Omega}}\big( 
    u^\varepsilon(x) - u(x)\big)
    ~\leq~ \overline{\mathbf{\Lambda}}\cdot \sqrt{\varepsilon},
\end{equation}
where the constant $\overline{\mathbf{\Lambda}}$ is explicitly computed as 
\bel{C-osc}
    \overline{\mathbf{\Lambda}} ~\doteq~
    \frac{2^{\alpha_p}}{\alpha_p} + \frac{1 + \left(n+2\right)\left(\left(\mathrm{osc}_\Omega f\right)^{1/p}\|f\|_{{\mathbf{Lip}}}\right)^{1/2}}{\lambda}~~.
\eeq
\end{proposition}

{\bf Proof.} The proof is divided into three main steps:

{\bf 1.} We first assume that $f\in \mathcal{C}_c(\Omega)\cap {\bf Lip}(\overline{\Omega})$  is nonnegative such that 
\bel{supp-f}
\mathrm{supp}(f)~\subseteq~\Omega_{\kappa}~~\mathrm{for~some}~ \kappa\in (0,2\delta_{\Omega}).
\eeq
In this case, we have  
\bel{va-u}
u^{\ve}(x)~\geq~u(x)~=~0\qquad\forall x\in U\doteq\Omega\backslash \overline{\Omega}_\kappa.
\eeq
Let $\tilde{u}^\varepsilon$  be the viscosity solution of \eqref{HJ-ve} with $f\equiv 0$ and $\Omega=U$. Since $u^{\ve}$ is also the subsolution of  \eqref{HJ-ve}, the comparison principle yields 
\[
\tilde{u}^{\ve}(x)~\geq~u^{\ve}(x)\qquad\forall x\in U.
\]
On the other hand, by  Lemma \ref{lem:ZeroData}, we get
\bel{ek1}
0~\leq~u^{\ve}(x)~\leq~\tilde{u}^{\ve}(x)~\leq~\left({1\over\lambda}+\frac{2^{\alpha_p}}{\alpha_p}\right)\cdot \ve^{1-\frac{1}{2}\alpha_p}, \qquad x\in U,
\eeq
and (\ref{va-u}) yields
\[
u^{\ve}(x)-u(x)~\leq~\left({1\over\lambda}+\frac{2^{\alpha_p}}{\alpha_p}\right)\cdot \ve^{1-\frac{1}{2}\alpha_p},\qquad x\in U.
\]

{\bf 2.} Next, we proceed with the standard doubling variable method by considering the  following auxiliary functional with a parameter $\gamma>0$ (taken to be large):
\begin{equation}\label{eq:auxiliary}
	\Phi^{\gamma}(x,y)~=~ u^\varepsilon(x) - u(y) -{\gamma\over 2}\cdot |x-y|^2 ,
	\qquad (x,y) \in \overline{\Omega}\times \overline{\Omega}.
\end{equation}
 Picking a global maximizer  $(x_\gamma, y_\gamma)\in \overline{\Omega}\times\overline{\Omega}$  of $\Phi^{\gamma}(x,y)$ over  $\overline{\Omega}\times\overline{\Omega}$, we have 
 \[
 u^\varepsilon(x_\gamma) - u(y_\gamma) - {\gamma\over 2}\cdot \big|x_{\gamma}-y_{\gamma}\big|^2~\geq~ u^\varepsilon(x_\gamma) - u(x_\gamma),
 \]
 and this yields 
 \[
 {\gamma\over 2}\cdot |x_{\gamma}-y_{\gamma}|^2~\leq~u(x_\gamma)-u(\gamma_\gamma)~\leq~\left(\mathrm{osc}_\Omega f\right)^{1/p}\cdot |x _{\gamma}-y_{\gamma}|.
 \]
Therefore
\begin{equation}\label{D-b}
    |x_{\gamma}-y_{\gamma}|~\leq~{2\left(\mathrm{osc}_\Omega f\right)^{1/p}\over \gamma},
\end{equation}
 where we invoke Lemma \ref{lem:BoundLipU} for the Lipschitz constant of $u$. Here two cases are considered:

\n$\bullet$ If $x_{\gamma}\in U$, then since $f=0$ in $U$, for every $x\in\Omega$ it holds that
\bel{epsalpha}
u^\varepsilon(x) - u(x)~=~\Phi(x,x)~\leq~ \Phi(x_\gamma, y_\gamma)~\leq~u^\varepsilon(x_\gamma) ~\leq~\left({1\over\lambda}+\frac{2^{\alpha_p}}{\alpha_p}\right)\cdot \ve^{1-\frac{1}{2}\alpha_p},\,
\eeq
where we use the fact that $u\geq 0$.

$\bullet$ Otherwise if $x_\gamma\in \overline{\Omega}_{\kappa}\subset\Omega$, then  since $x_\gamma$ is a  maximizer of $\Phi^{\gamma}(x,y_{\gamma})$ over $\Omega$, it follows that
\[
|Du^{\ve}(x_\gamma)|~=~\gamma|x_{\gamma}-y_{\gamma}|,\qquad \Delta u^{\ve}(x_{\gamma})~=~\Delta_{x}\left(\Phi(x_{\gamma},y_{\gamma})+{\gamma\over 2}\cdot |x-y_{\gamma}|^2\right)~\leq~\gamma n.
\]
Thus, (\ref{HJ-ve}) implies that
\[
\begin{split}
u^{\ve}(x_{\gamma})&~=~{1\over \lambda}\cdot\Big(f(x_{\gamma})+\ve \Delta u^{\ve}(x_{\gamma})-|Du^{\ve}(x_\gamma)|^p\Big)~\leq~{1\over\lambda}\cdot\Big(f(x_{\gamma})+\ve \gamma n-\gamma^{p} |x_{\gamma}-y_{\gamma}|^p\Big).
\end{split}
\]
On the other hand, since $y_{\gamma}$ is a maximizer of $y\mapsto \Phi^{\gamma}(x_{\gamma},y)=\ds \left(u^{\ve}(x_{\gamma})-{\gamma\over 2}\cdot |x_{\gamma}-y|^2 \right)-u(y)$,  we can apply the supersolution test for $u$ in (\ref{HJ-S}) to obtain 
\[
u(y_{\gamma})~\geq~{1\over \lambda}\cdot \Big(f(y_{\gamma})-\gamma^{p} |x_{\gamma}-y_{\gamma}|^p\Big).
\]
Hence, by (\ref{D-b}) we derive
\[
u^{\ve}(x_{\gamma})-u(y_{\gamma})~\leq~{1\over \lambda}\cdot\Big(f(x_{\gamma})-f(y_{\gamma})+\ve \gamma n\Big)~\leq~{\|f\|_{{\mathbf{Lip}}}\cdot |x_{\gamma}-y_{\gamma}|+\ve n \gamma \over \lambda},
\]
which yields 
\[
\begin{split}
u^\varepsilon(x) - u(x)~=~\Phi^{\gamma}(x,x)&~\leq~\Phi^{\gamma}(x_{\gamma},y_{\gamma})
~\leq~{\|f\|_{{\mathbf{Lip}}}\cdot |x_{\gamma}-y_{\gamma}|+\ve n \gamma \over \lambda}-{\gamma\cdot |x_\gamma-y_{\gamma}|^2\over 2}\\
&~\leq~ \frac{1}{\lambda} \left( \ve n \gamma + \frac{2\left(\mathrm{osc}_\Omega f\right)^{1/p}\|f\|_{{\mathbf{Lip}}}}{\gamma}\right) \quad\forall x\in \Omega.
\end{split}
\]
Finally, if $\Vert f\Vert_{\mathbf{Lip}}\neq 0$, choosing 
\begin{equation*}
	\gamma = \left(\frac{(\mathrm{osc}_\Omega f)^{1/p}\Vert f\Vert_{\mathbf{Lip}}}{\varepsilon}\right)^{1/2},
\end{equation*}
we deduce that
\bel{eps12}
	u^\varepsilon(x) - u(x)~\leq~ \left(\frac{n+2}{\lambda} \right)\cdot \left( \left(\mathrm{osc}_\Omega f\right)^{1/p} \|f\|_{{\mathbf{Lip}}}\right)^{1/2}\cdot\varepsilon^{1/2}, \qquad\forall x\in \Omega.
\eeq
Then equations \eqref{epsalpha} and \eqref{eps12} yield \eqref{eq:THMEstimateUepsConstantfUpper}, since $\ds1-\frac{\alpha_p}{2} > \frac{1}{2}$ for every $p>2$.
\v

{\bf 3.}  Finally, to remove the assumption of $f$ having a compact support, for $0<\kappa<\delta_{\Omega}$ we consider a modification of $f$ defined by 
 \bel{f-k}
 f_{\kappa}(x)~=~\eta_{\kappa}(d_{\partial\Omega}(x))\cdot f(x)\qquad\forall x\in \overline{\Omega}
 \eeq
 where $\eta_{\kappa}:[0,\infty)\to[0,1]$ is a cut-off function such that 
\begin{equation}\label{eq:cutoff}
	0~\leq~\eta'_k~\leq~{4\over \kappa}~~,\qquad\quad
 \eta_{k}(s)~=~\begin{cases}
1&\qquad\mathrm{if}\qquad s\geq \kappa\\[2mm]
0&\qquad \mathrm{if}\qquad 0\leq s\leq \frac{\kappa}{2}
 \end{cases}~~.
\end{equation}
Since $f=0$ on $\partial\Omega$, we have 
\bel{fkLip}
\begin{cases}
\mathrm{supp}(f_k)\subseteq \Omega_{\kappa/2},\qquad\| f_{\kappa}\|_{{\mathbf{Lip}}}~\leq~\| f\|_{{\mathbf{Lip}}} + 2\kappa\\[3mm]
\ds\sup_{x\in \overline{\Omega}} |f(x)-f_{\kappa}(x)|~\leq~\sup_{x\in \overline{\Omega}\backslash\Omega_{\kappa}} |f(x)|~\leq~\|f\|_{{\mathbf{Lip}}}\cdot \kappa
\end{cases}~~.
 \eeq
Let \( u^\varepsilon_\kappa \in \mathcal{C}^2(\Omega) \cap \mathcal{C}^{0,\alpha}(\overline{\Omega}) \) and $u_\kappa\in \mathrm{Lip}(\overline{\Omega})$ be the respective solutions to \eqref{HJ-ve} and \eqref{HJ-S} with \( f \) replaced by \( f_\kappa \). By the comparison principle in Proposition \ref{prop:Comparison} for the H\"older subsolution and supersolution of \eqref{HJ-ve}, and the standard comparison principle for \eqref{HJ-S}, we obtain
\begin{equation}\label{eq:EstimatesKappa}
	 0~\leq~  u^\varepsilon(x) - u^\varepsilon_\kappa(x), u(x) - u_\kappa(x)~\leq~ \frac{\|f\|_{{\mathbf{Lip}}}}{\lambda}\cdot\kappa,\qquad x\in \overline{\Omega}. 
\end{equation}
It also follows that 
\begin{equation}\label{eq:EstimatesKappaUEps}
	u^\varepsilon_\kappa(x) - u_\kappa(x)~\leq~ C(\mathrm{osc}_\Omega f_\kappa,p,n)\cdot \varepsilon^{\frac{1}{2}}, \qquad x\in \overline{\Omega}~.
\end{equation}
From \eqref{eq:EstimatesKappa} and \eqref{eq:EstimatesKappaUEps} we get
\begin{align*}
	u^\varepsilon(x) - u(x) 
	&~=~ \Big( u^\varepsilon(x) - u^\varepsilon_\kappa(x) \Big) + 
	\Big(u^\varepsilon_\kappa(x) - u_\kappa(x) \Big) + \Big( u_\kappa(x) - u(x)\Big) \\
	&~\leq~  \frac{\|f\|_{{\mathbf{Lip}}}}{\lambda}\cdot\kappa+ \overline{\mathbf{\Lambda}}  \cdot \varepsilon^{\frac{1}{2}}, \qquad x\in \overline{\Omega}.
\end{align*}  
Taking $\kappa\to 0^+$, we obtain \eqref{eq:THMEstimateUepsConstantfUpper}, since $\mathrm{osc}_\Omega f_\kappa \to \mathrm{osc}_\Omega f$ as $\kappa\to 0^+$ due to \eqref{fkLip}.
\endproof

\subsection{Upper bound of $u^{\ve}-u$ for semiconcave data}\label{con-se}
To begin the proof of Theorem \ref{semi-im}, we first establish a simple lemma on the semiconcavity of the solution in both the subquadratic and superquadratic cases, assuming that the data $f$ is nonnegative and semiconcave. This type of result was previously studied for $1 < p\leq 2$ in  \cite{ han_global_2025, han_remarks_2022}. We refer to \cite{{bardi_optimal_1997}, tran_hamilton-jacobi_2021} for more details on the optimal control formula of solutions to Hamilton-Jacobi equations.  

\begin{lemma}\label{semi-c} For every $p>1$, assume that  $f\in \mathcal{C}_c(\overline{\Omega})$ is nonnegative and semiconcave with a semiconcavity constant $\mathbf{c}_f$. Then the solution $u$ of (\ref{HJ-S}) is semiconcave with a semiconcavity constant $\mathbf{c}_f/\lambda$.
\end{lemma}
{\bf Proof.} Fix any $x\in \Omega$. By (\ref{L}) and (\ref{O-C}), we have
\bel{u-x}
u(x)~=~\inf \left\lbrace 
	\int_{-\infty}^{0}e^{\lambda s} 
	\big(c_p |\dot{\eta}(s)|^q+f(\eta(s))\big)~ds: 
    \eta\in \mathcal{AC}((-\infty,0];\overline{\Omega}), \eta(0)=x
    \right\rbrace.
\eeq
Let $\eta\in \mathcal{AC}((-\infty,0];\overline{\Omega})$ be a minimizer to \eqref{u-x} with $\eta(0) = x$.  As $f\in  \mathcal{C}_c(\overline{\Omega})$ is  nonnegative, the optimality condition yields
\bel{cs1}
	\eta(s)~\in~\overline{\mathrm{supp}(f)}\qquad\forall s\leq 0.
\eeq
Indeed, since $u = 0$ whenever $f = 0$, only points $x \in \mathrm{supp}(f)$ serve as relevant starting positions. For these points, the optimal trajectory in \eqref{u-x} stays within the closure of $\mathrm{supp}(f)$. Once it reaches a region where $f = 0$, the trajectory may simply remain stationary.
For every $\tau<0$ and $h\in \R^n$,  consider the curves
\begin{align*}
    \eta^{\pm}_{\tau,h}(s)~ = ~
    \begin{cases}
     \ds   \eta(s) + \left(1\pm\frac{s}{\tau}\right)\cdot h, &\;\;\;\; \tau\leq s\leq 0 \\[3mm]
        \eta(s), & -\infty\leq s < \tau
    \end{cases}~~.
\end{align*}
By (\ref{cs1}), for  sufficiently small $|h|>0$, we have  that $\eta^\pm_{\tau,h} \in \mathcal{AC}((-\infty,0]; \Omega)$ and $\eta^\pm_{\tau,h}(0) = x \pm h$ for all $\tau<0$ . Thus, by the optimality condition and the semiconcavity property of $f$, we estimate 
\[
\begin{split}
u(x+h)+u(x-h)-2u(x)&~\leq~\int_{\tau}^0 c_pe^{\lambda s}\cdot \left(|\dot{\eta}_{\tau,h}^{+}|^q+|\dot{\eta}_{\tau,h}^{-}|^q-2|\dot{\eta}|^q\right)ds\\
&\qquad\qquad\qquad +\int_{\tau}^0e^{\lambda s}\cdot \big[f(\eta_{\tau,h}^+)+f(\eta_{\tau,h}^-)-2f(\eta)\big]ds\\
&~\leq~\O(1)\cdot {|h|^q\over |\tau|^{q-1}}+{\bf c}_f|h|^2\cdot \int_{\tau}^0e^{\lambda s} \left(1+{s\over\tau}\right)^2ds.
\end{split}
\]
Taking $\tau\to-\infty$, we obtain 
\[
u(x+h)+u(x-h)-2u(x)~\leq~{{\bf c}_f\over \lambda}\cdot |h|^2,
\]
which completes the proof. 
\endproof
\v

{\bf Proof of Theorem \ref{semi-im}.} Assume that $f$ is semiconcave with a  semiconcavity constant $\mathbf{c}_f$ and $\mathrm{supp}(f)\subset\Omega_{\delta_f}$ for some $\delta_f,\mathbf{c}_f>0$. As in the first step of the proof of Proposition \ref{prop:upper}, it holds
\bel{ma1}
u^{\ve}(x)-u(x)~\leq~\left({1\over\lambda}+\frac{2^{\alpha_p}}{\alpha_p}\right)\cdot \ve^{1-\frac{1}{2}\alpha_p},\qquad x\in U_f\doteq\Omega\backslash\delta_f.
\eeq
 Let $x_{\max}\in \overline{\Omega}$ be  such that 
\bel{dq}
   u^\varepsilon(x_{\max}) - u(x_{\max})~=~ \max_{x\in \overline{\Omega}} \big\{u^\varepsilon(x)- u(x)\big\}. 
\eeq
If $x_{\max}\in \overline{\Omega}_{\delta_f}$, then using \( u^\varepsilon \) as a test function in the subsolution condition for \( u \) at \( x_{\max} \) and applying the equation for \( u^\varepsilon \) at \( x_{\max} \in \Omega \), we get
\[
\begin{cases}
\lambda u(x_{\max})+|D u^{\ve}(x_{\max})|^{p}-f(x_{\max})~\geq~0\\[2mm]
\lambda u^{\ve}(x_{\max})+|D u^{\ve}(x_{\max})|^{p}-f(x_{\max})~=~\ve\cdot \Delta u^{\ve}(x_{\max}),
\end{cases}
\]
which implies
\[
 u^{\ve}(x_{\max})- u(x_{\max})~\leq~{\ve\over \lambda}\cdot  \Delta u^{\ve}(x_{\max}).
\]
Recalling Lemma \ref{semi-c}, we have that  $u$  is semiconcave with a  semiconcavity constant $\ds{\mathbf{c}_f\over \lambda}$ and 
\[
 u^{\ve}(x_{\max})- u(x_{\max})~\leq~{\ve\over \lambda}\cdot  \Delta u^{\ve}(x_{\max})~\leq~{\ve\over \lambda}\cdot{\mathbf{c}_f\over \lambda}~=~{{\bf c}_f\over \lambda^2}\cdot \ve~.
\]
Hence from (\ref{ma1})-(\ref{dq}), we have 
\[
\max_{x\in \overline{\Omega}} \big\{u^\varepsilon(x)- u(x)\big\}~\leq~\max\left\{\left({1\over\lambda}+\frac{2^{\alpha_p}}{\alpha_p}\right)\cdot \ve^{1-\frac{1}{2}\alpha_p},{{\bf c}_f\over \lambda^2}\cdot \ve\right\},
\]
and (\ref{ke}) holds for every $\ve>0$ sufficiently small.
\endproof

\begin{corollary}\label{sf}
For every $p>2$, assume that $\Omega\subset\R^n$ is an open and bounded subset of $\R^n$ with $\mathcal{C}^2$ boundary, $f\in \mathcal{C}^{2}(\overline{\Omega})$ is nonnegative, and 
\bel{Bdc}
f(x)~=~0,\qquad Df(x)~=~0~~~\mathrm{on}~\partial\Omega~.
\eeq
Then for every $\ve>0$ sufficiently small, it holds that
\[
\max_{x\in\overline{\Omega}}\{u^{\ve}(x)-u(x)\}~\leq~\left({1\over\lambda}+\frac{2^{\alpha_p}}{\alpha_p}\right)\cdot \ve^{1-\frac{\alpha_p}{2}}~.
\]
\end{corollary}
{\bf Proof.}  By  \cite[Lemma 18]{han_remarks_2022}, there exists a sequence of nonnegative functions $f_k\in C_c(\Omega)$ such that $f_k$ is semiconcave with a  semiconcavity constant $K_f>0$ and 
\[
\mathrm{supp}(f_k)~\subset~\Omega_{1/k},\qquad \lim_{k\to\infty}\|f-f_k\|_{\infty}~=~0.
\]
Let \( u^\varepsilon_k \in \mathcal{C}^2(\Omega) \cap \mathcal{C}^{0,\alpha}(\overline{\Omega}) \) and \( u_k \in \mathbf{Lip}(\overline{\Omega}) \) be the solutions to \eqref{HJ-ve} and \eqref{HJ-S}, respectively, with \( f \) replaced by \( f_k \). By comparison principle, we have 
\bel{l-e1}
\lim_{k\to\infty} \|u_k-u\|_{\infty}~=~\lim_{k\to\infty} \big\|u^{\ve}_k-u^\ve\big\|_{\infty}~=~0.
\eeq
Set $U_k\doteq\Omega\backslash\overline{\Omega}_{1/k}$. As in the proof of  Theorem \ref{semi-im}, for every $\ve>0$ sufficiently small, it holds that
\[
\max_{x\in \overline{\Omega}} \big\{u^\varepsilon_k(x)- u_k(x)\big\}~\leq~\max\left\{\left({1\over\lambda}+\frac{2^{\alpha_p}}{\alpha_p}\right)\cdot \ve^{1-\frac{1}{2}\alpha_p},{K_f\over \lambda^2}\cdot \ve\right\}~=~\left({1\over\lambda}+\frac{2^{\alpha_p}}{\alpha_p}\right)\cdot \ve^{1-\frac{1}{2}\alpha_p}~.
\]
Taking $k\to\infty$ and using (\ref{l-e1}), we achieve
\[
\max_{x\in\overline{\Omega}}\{u^{\ve}(x)-u(x)\}~\leq~\limsup_{k\to\infty}\max_{x\in \overline{\Omega}} \big\{u^\varepsilon_k(x)- u_k(x)\big\}~\leq~\left({1\over\lambda}+\frac{2^{\alpha_p}}{\alpha_p}\right)\cdot \ve^{1-\frac{1}{2}\alpha_p},
\]
which completes the proof.
\endproof
\v

{\bf Acknowledgement.} This research by Khai T. Nguyen was partially supported by National Science Foundation grant DMS-2154201.


\begin{thebibliography}{10}

\bibitem{armstrong_viscosity_2015}
{\sc Armstrong, S.~N., and Tran, H.~V.}
\newblock Viscosity solutions of general viscous {Hamilton}–{Jacobi}
  equations.
\newblock {\em Mathematische Annalen 361}, 3 (2015), 647--687.

\bibitem{bardi_optimal_1997}
{\sc Bardi, M., and Capuzzo-Dolcetta, I.}
\newblock {\em Optimal {Control} and {Viscosity} {Solutions} of
  {Hamilton}–{Jacobi}–{Bellman} {Equations}}.
\newblock Modern {Birkhäuser} {Classics}. Birkhäuser Basel, 1997.

\bibitem{barles_short_2010}
{\sc Barles, G.}
\newblock A short proof of the $c^{0,\alpha}$-regularity of viscosity
  subsolutions for superquadratic viscous {Hamilton}–{Jacobi} equations and
  applications.
\newblock {\em Nonlinear Analysis: Theory, Methods \& Applications 73}, 1 (2010), 31--47.

\bibitem{barles_generalized_2004}
{\sc Barles, G., and Da~Lio, F.}
\newblock On the generalized {Dirichlet} problem for viscous
  {Hamilton}–{Jacobi} equations.
\newblock {\em Journal de Mathématiques Pures et Appliquées 83}, 1 (2004), 53--75.

\bibitem{cannarsa_semiconcave_2004}
{\sc Cannarsa, P., and Sinestrari, C.}
\newblock {\em Semiconcave {Functions}, {Hamilton}—{Jacobi} {Equations}, and
  {Optimal} {Control}}.
\newblock Birkhäuser, Boston, MA, 2004.

\bibitem{capuzzo-dolcetta_holder_2010}
{\sc Capuzzo-Dolcetta, I., Leoni, F., and Porretta, A.}
\newblock Hölder estimates for degenerate elliptic equations with coercive
  {Hamiltonians}.
\newblock {\em Transactions of the American Mathematical Society 362}, 9 (2010), 4511--4536.

\bibitem{capuzzo-dolcetta_hamilton-jacobi_1990}
{\sc Capuzzo-Dolcetta, I., and Lions, P.-L.}
\newblock Hamilton-{Jacobi} {Equations} with {State} {Constraints}.
\newblock {\em Transactions of the American Mathematical Society 318}, 2
  (1990), 643--683.


\bibitem{chaintron_optimal_2025}
{\sc Chaintron, L.-P., and Daudin, S.}
\newblock Optimal rate of convergence in the vanishing viscosity for quadratic
  {Hamilton}-{Jacobi} equations, 2025.
\newblock arXiv:2502.09103 [math].

\bibitem{cirant_convergence_2025}
{\sc Cirant, M., and Goffi, A.}
\newblock Convergence rates for the vanishing viscosity approximation of
  {Hamilton}-{Jacobi} equations: the convex case, 2025.
\newblock arXiv:2502.15495 [math].

\bibitem{crandall_users_1992}
{\sc Crandall, M.~G., Ishii, H., and Lions, P.-L.}
\newblock User's guide to viscosity solutions of second order partial
  differential equations.
\newblock {\em Bull. Amer. Math. Soc. (N.S.) 27} (1992), 1--67.

\bibitem{crandall_two_1984}
{\sc Crandall, M.~G., and Lions, P.~L.}
\newblock Two {Approximations} of {Solutions} of {Hamilton}-{Jacobi}
  {Equations}.
\newblock {\em Mathematics of Computation 43}, 167 (1984), 1--19.


\bibitem{fleming_convergence_1961}
{\sc Fleming, W.~H.}
\newblock The convergence problem for differential games.
\newblock {\em Journal of Mathematical Analysis and Applications 3}, 1 (1961), 102--116.

\bibitem{gilbarg_elliptic_2001}
{\sc Gilbarg, D., and Trudinger, N.~S.}
\newblock {\em Elliptic {Partial} {Differential} {Equations} of {Second}
  {Order}}, 2~ed.
\newblock Classics in {Mathematics}. Springer-Verlag, Berlin Heidelberg, 2001.

\bibitem{han_global_2025}
{\sc Han, Y.}
\newblock Global semiconcavity of solutions to first-order {Hamilton}-{Jacobi}
  equations with state constraints.
\newblock {\em Advances in Continuous and Discrete Models}, 106 (2025),
  1--36.

\bibitem{han_quantitative_2025}
{\sc Han, Y., Jing, W., Mitake, H., and Tran, H.~V.}
\newblock Quantitative {Homogenization} of {State}-{Constraint}
  {Hamilton}–{Jacobi} {Equations} on {Perforated} {Domains} and
  {Applications}.
\newblock {\em Archive for Rational Mechanics and Analysis 249}, 18 (2025),
  1--53.

\bibitem{han_remarks_2022}
{\sc Han, Y., and Tu, S. N.~T.}
\newblock Remarks on the {Vanishing} {Viscosity} {Process} of
  {State}-{Constraint} {Hamilton}–{Jacobi} {Equations}.
\newblock {\em Applied Mathematics \& Optimization 86}, 3 (2022), 1--42.

\bibitem{kim_state-constraint_2020}
{\sc Kim, Y., Tran, H.~V., and Tu, S.~N.}
\newblock State-{Constraint} {Static} {Hamilton}--{Jacobi} {Equations} in
  {Nested} {Domains}.
\newblock {\em SIAM Journal on Mathematical Analysis 52}, 5 (2020),
  4161--4184.


\bibitem{lasry_nonlinear_1989}
{\sc Lasry, J.~M., and Lions, P.-L.}
\newblock Nonlinear {Elliptic} {Equations} with {Singular} {Boundary}
  {Conditions} and {Stochastic} {Control} with {State} {Constraints}. {I}.
  {The} {Model} {Problem}.
\newblock {\em Mathematische Annalen 283}, 4 (1989), 583--630.


\bibitem{porretta_local_2004}
{\sc Porretta, A.}
\newblock Local estimates and large solutions for some elliptic equations with
  absorption.
\newblock {\em Advances in Differential Equations 9}, 3-4 (2004),
  329--351.

\bibitem{alessio_asymptotic_2006}
{\sc {Porretta}, A., and {Véron}, L.}
\newblock Asymptotic {Behaviour} of the {Gradient} of {Large} {Solutions} to
  {Some} {Nonlinear} {Elliptic} {Equations}.
\newblock {\em Advanced Nonlinear Studies 6}, 3 (2006), 351--378.

\bibitem{qian_optimal_2024}
{\sc Qian, J., Sprekeler, T., Tran, H.~V., and Yu, Y.}
\newblock Optimal {Rate} of {Convergence} in {Periodic} {Homogenization} of
  {Viscous} {Hamilton}-{Jacobi} {Equations}.
\newblock {\em Multiscale Modeling \& Simulation 22}, 4 (2024),
  1558--1584.


\bibitem{soner_optimal_1986}
{\sc Soner, H.~M.}
\newblock Optimal {Control} with {State}-{Space} {Constraint} {I}.
\newblock {\em SIAM Journal on Control and Optimization 24}, 3 (1986),
  552--561.


\bibitem{tabet_tchamba_large_2010}
{\sc Tabet~Tchamba, T.}
\newblock Large time behavior of solutions of viscous {Hamilton}–{Jacobi}
  equations with superquadratic {Hamiltonian}.
\newblock {\em Asymptotic Analysis 66}, 3-4 (2010), 161--186.


\bibitem{tran_adjoint_2011}
{\sc Tran, H.~V.}
\newblock Adjoint methods for static {Hamilton}–{Jacobi} equations.
\newblock {\em Calculus of Variations and Partial Differential Equations 41}, 3
  (2011), 301--319.

\bibitem{tran_hamilton-jacobi_2021}
{\sc Tran, H.~V.}
\newblock {\em Hamilton-{Jacobi} {Equations}: {Theory} and {Applications}}.
\newblock American Mathematical Society, 2021.
\newblock Google-Books-ID: H2tMzgEACAAJ.

\end{thebibliography}
\end{document}